\documentclass[reqno]{amsart}
\usepackage[utf8]{inputenc}

\usepackage[utf8]{inputenc}

\usepackage{sky}
\usepackage{graphicx}
\usepackage{color}
\usepackage{caption}
\usepackage{subcaption}

\newcommand{\vertiii}[1]{{\left\vert\kern-0.25ex\left\vert\kern-0.25ex\left\vert #1 
\right\vert\kern-0.25ex\right\vert\kern-0.25ex\right\vert}}

\title{$\mrm{U}(N)$ lattice Yang--Mills in the 't~Hooft regime}

\author{Ron Nissim}
\address{Ron Nissim, Department of Mathematics, Massachusetts Institute of Technology, Cambridge, MA 02139}
\email{rnissim@mit.edu}

\date{}

\begin{document}

\begin{abstract}
    We establish a mass gap, prove the existence of a unique infinite volume limit, and give a new proof of the large $N$ limit for $\mathrm{U}(N)$ lattice Yang--Mills theory in the 't~Hooft regime. These results were previously obtained for $\mathrm{SU}(N)$ and $\mathrm{SO}(N)$ lattice Yang–Mills theories as applications of the mixing of the associated Langevin dynamics, which is verified via the Bakry--Émery criterion \cite{shen2023stochastic}. For $\mathrm{U}(N)$, however, this approach fails because its Ricci curvature is not uniformly positive, and as a result the Bakry--Émery condition cannot be easily verified. To overcome this obstacle, we recast the $\mathrm{U}(N)$ theory as a random-environment $\mathrm{SU}(N)$ model, where the randomness arises from a $\mathrm{U}(1)$ field, and combine cluster-expansion and Langevin-dynamics techniques to analyze the resulting $\mathrm{U}(1)\times\mathrm{SU}(N)$ model.
\end{abstract}

\maketitle

\tableofcontents

\section{Introduction}

Euclidean quantum Yang--Mills theory is the mathematical framework for the standard model of particle physics. Unfortunately the model is not well-defined directly in the continuum. In order to initiate a rigorous mathematical treatment of the subject, \cite{Wilson1974} introduced a discretized lattice version of the model referred to as lattice Yang--Mills theory. The data of the model is a compact Lie group $G$ known as the \textit{gauge group}, a finite subgraph $\Lambda \subset \Z^d$, and an \textit{inverse coupling strength} $\beta>0$ (We always use the 't Hooft scaling for $\beta$ which will be defined later). Given this data, there is an associated lattice Yang-Mills probability measure $\mu_{G,\Lambda,\beta}$ on the space of assignments, $Q:E_{\Lambda}^+ \to G$, of group elements to each (positively oriented) edge. We refer the reader to \cite{chatterjee2016a} for a more complete survey of the model and the important open problems in the subject.

We say that the Yang--Mills measure obeys a \textit{mass gap} if the covariance between two observables, which only depend on edge variables separated by a distance $R$, decays exponentially in $R$. Establishing a mass gap in four dimensions is one of the major open problems in Yang--Mills theory. In particular, a lattice based approach to solving the Yang--Mills millennium problem \cite{Jaffe2006a} would require proving a mass gap for arbitrarily large inverse coupling parameter $\beta$ for a lattice Yang-Mills theory with a non-Abelian gauge group $G$. One important application of a mass gap, is the construction of a unique infinite volume limit for the measure $\mu_{G,\Lambda,\beta}$ as $\Lambda \to \Z^d$.

While it is suspected that Yang-Mills theories, such as those with groups $G \in \{\mrm{SU}(N),\mrm{SO}(N),\mrm{U}(N)\}$ for $N \geq 3$ exhibit a mass gap at all $\beta>0$ for $d\in \{3,4\}$, the existing results always involve a strong coupling (small $\beta$) assumption. In the classic work of Osterwalder and Seiler \cite{osterwalderseiler1978}, they apply a cluster expansion argument to establish that every lattice Yang-Mills theory, regardless of gauge group $G$ and dimension $d$, has a mass gap for sufficiently small $\beta>0$. If one tracks their arguments for the groups $G=\{\mrm{SU}(N),\mrm{SO}(N),\mrm{U}(N)\}$, their condition reduces to $\beta < \frac{c_d}{N^2}$ for some small dimension dependent constant $c_d$. In contrast, for Abelian Lattice Yang-Mills when $G=\mrm{U}(1)$, it is known that in 3D there is a mass gap for all $\beta>0$ \cite{gopfertmack1981}, while for 4D there is a known phase transition so that mass gap does not hold for sufficiently large $\beta$ \cite{frohlich1982massless}.

Recently a new perspective on lattice Yang--Mills at strong coupling has emerged based on an associated dynamics. In the work of \cite{shen2023stochastic} the authors use this dynamical approach to establish a mass gap and unique infinite volume for $G \in \{\mrm{SU}(N),\mrm{SO}(N)\}$ in the strong coupling regime when $\beta < c_d$ for some explicit dimensional constant $c_d$, significantly improving on the strong coupling regime of \cite{osterwalderseiler1978}. Under the scaling used throughout the paper, this strong coupling regime when $\beta < c_d$ is known as the 't Hooft regime. In follow up work the same authors extend their results to the Yang--Mills--Higgs model \cite{SZZ2024}, and finally in a recent work the dynamical approach was also used to establish \textit{area law} \cite{CNS2025b} in the same strong coupling regime for $G \in \{\mrm{SU}(N),\mrm{SO}(2N),\mrm{U}(N)\}$. 

Despite the achievements of the dynamical approach, it is limited in the groups for which it directly applies. In particular while the arguments of \cite{shen2023stochastic} apply for $G \in \{\mrm{SU}(N),\mrm{SO}(N)\}$, they do not carry over to $G=\mrm{U}(N)$. In this paper we express the $\mrm{U}(N)$ as an $\mrm{U}(1) \times \mrm{SU}(N)$ model and use a combination of ideas from the cluster expansion of \cite{osterwalderseiler1978} and the dynamical approach to prove the mass gap and unique infinite volume for $\mrm{U}(N)$ Lattice Yang--Mills theory in the 't Hooft regime. We state this precisely in the following theorem, referring the reader to Section \ref{subsection:Lattice-YM-Prelim} for the precise notation.

\begin{theorem}[Mass gap and infinite volume limit]\label{thm:U-N-mass-gap}
Let $d \geq 2$, $N \geq 2$. Then for some fixed $\tilde{\beta}=\tilde{\beta}(d)$, and all $\beta < \tilde{\beta}$, there exists a probability measure $\mu_{\mrm{U}(N),\beta}$ on $\mrm{U}(N)^{E_{\Lambda_{\infty}}^+}$ (with $\Lambda_{\infty}=\Z^d$) such that
\begin{equ}
    \mu_{\mrm{U}(N),\Lambda_L,\beta} \to \mu_{\mrm{U}(N),\beta}
\end{equ}
weakly as $L \to \infty$. Additionally, for any smooth local observables $f, g$ with $\Lambda_f, \Lambda_g \sse \Lambda$ (see Notation \ref{notation:local-observables}), there exist constants $C_1=C_1(N,d,|\Lambda_f|,|\Lambda_g|)$ and $C_2=C_2(N,d)$ such that,
    \begin{equs}\label{eq:main-mass-gap}
        \big| \mathrm{Cov}_{\mu_{\mrm{U}(N),\beta}}(f,g)\big| \leq C_1 \big(\|f\|_{L^{\infty}}+\vertiii{f}_{\infty}\big)\big(\|g\|_{L^{\infty}}+\vertiii{g}_{\infty}\big) e^{-C_2 d(\Lambda_f,\Lambda_g)},
    \end{equs}
where $d(\Lambda_f, \Lambda_g)$ is the graph distance between the subsets $\Lambda_f, \Lambda_g$ within $\Lambda$.
\end{theorem}
\begin{remark}
    In order to prove Theorem \ref{thm:U-N-mass-gap}, we will first prove the mass gap statement \eqref{eq:main-mass-gap} on the finite lattice $\Lambda_L$ with constants independent of $L$. Then later after establishing the infinite volume limit, we can simply take the limit on both sides of the finite $L$ version of \eqref{eq:main-mass-gap} to recover \eqref{eq:main-mass-gap}.
\end{remark}

\begin{remark}
In fact, our proof of the infinite volume limit result works equally well for $G \in \{\mrm{SO}(N),\mrm{SU}(N)\}$ in the 't Hooft regime. Our argument does not rely on a sophisticated coupling argument like the one in \cite{shen2023stochastic}. However, our proof only shows that $\mu_{N,\Lambda_L,\beta}$ converges weakly, and in particular it does not show that the infinite volume dynamics has a unique invariant measure.
\end{remark}

\begin{remark}
    The proof of the mass gap given in Section \ref{section:Mass-Gap}  extends to the more general mixed-temperature model where each plaquette $p$ is assigned a possibly different inverse coupling parameter $\beta_p$ as long as $\sup_p \beta_p < \tilde{\beta}$. We state this result more precisely in Lemma \ref{lm:Mixed-Temp-Mass-Gap} and use it to establish the unique infinite volume limit part of Theorem \ref{thm:U-N-mass-gap}.
\end{remark}

There is additional significance to the 't Hooft regime. In particular with the 't Hooft scaling, for small fixed $\beta$, after sending $N \to \infty$, it has long been known that the expectations and correlations for the main observables of interest, \textit{Wilson loop observables}, simplify \cite{Hooft1973}. Moreover one can rigorously construct a limiting theory as $N\to \infty$ in the sense of assigning values to the expectations of Wilson loop observables. These expectations can be expressed by string theoretic sums, or sums over surfaces \cite{Chatterjee2019a,jafarov2016,BCSK2024,BCSK25}. One particular property of the large $N$ limit known as Wilson loop factorization asserts that in the $N\to \infty$ limit, the Wilson loop observables become completely deterministic, and hence uncorrelated with each other. While this fact can be established via the string theoretic or surface sum approach, in \cite{shen2023stochastic} it was observed that a short proof of this property could also be given as a consequence of the dynamical approach. Once again their proof only applies for $G \in \{\mrm{SU}(N),\mrm{SO}(N)\}$. We are able to extend their Large $N$ limit result to $\mrm{U}(N)$ by viewing $\mrm{U}(N)$ Lattice Yang--Mills theory as an $\mrm{U}(1) \times \mrm{SU}(N)$ model, and proving a Poincaré inequality for the marginal law corresponding to the $\mrm{SU}(N)$ field. We state the precise theorem below, once again referring the reader to Section \ref{subsection:Lattice-YM-Prelim} for the precise notation.

\begin{theorem}[Large $N$ Limit for $\mrm{U}(N)$]\label{thm:main-large-N-Limit}
    Suppose $\beta < \tilde{\beta}(d)$. Then for any loop $\ell$ in $\Lambda_{\infty}=\Z^d$,
    \begin{equs}
        W_{\ell}-\langle W_{\ell}\rangle_{\mu_{\mrm{U}(N),\beta}} \to 0
    \end{equs}
    in probability, and as a consequence for a collection of loops $\ell_1,\dots, \ell_n$, we have the following Wilson loop factorization property,
    \begin{equs}
        \langle W_{\ell_1}W_{\ell_2}\cdots W_{\ell_n}\rangle_{\mu_{\mrm{U}(N),\beta}}- \langle W_{\ell_1}\rangle_{\mu_{\mrm{U}(N),\beta}}\langle W_{\ell_2}\rangle_{\mu_{\mrm{U}(N),\beta}}\cdots \langle W_{\ell_n}\rangle_{\mu_{\mrm{U}(N),\beta}} \to 0
    \end{equs}
\end{theorem}
\begin{remark}
    When comparing to \cite[Corollary 1.5]{shen2023stochastic}, note that their definition of a Wilson loop observable differs from ours by a multiplicative factor $N$.
\end{remark}

We now elaborate on why the techniques of \cite{shen2023stochastic} don't apply for $G=\mrm{U}(N)$. The key to their dynamical approach is to show that the Langevin dynamics they study mixes exponentially fast. This rapid mixing property is typically proven by verifying a convexity condition known as the Bakry-Émery condition \cite{GuionnetLogSobolev}. The verification of the Bakry-Émery condition in \cite{shen2023stochastic} critically relies on the good positive lower bound on the Ricci curvature of $\mrm{SU}(N)$ and $\mrm{SO}(N)$. This positive lower bound on the Ricci curvature fails for $\mrm{U}(N)$. As one doesn't expect the Bakry-Émery condition to hold, to treat the $\mrm{U}(N)$ case one needs to either verify that the Langevin dynamics mixes through an alternate approach, or modify the dynamics itself. We use the latter approach.

We now explain the key ideas for the proof of the mass gap statement of Theorem \ref{thm:U-N-mass-gap}. Our arguments rely on treating $\mrm{U}(N)$ lattice Yang--Mills theory as a sort of random environment $\mrm{SU}(N)$ lattice Yang--Mills theory, where the random environment is itself some $\mrm{U}(1)$ random field. We also remark that in the treatment of the dynamical treatment of the Yang-Mills-Higgs model \cite{SZZ2024}, the authors consider a similar decomposition of the measure. With our decomposition in mind, the argument to establish a mass gap is broken down into three steps. By first conditioning on the $\mrm{SU}(N)$ part of the measure, we are able to establish a conditional mass gap result using a cluster expansion argument similar to \cite{osterwalderseiler1978}. One difference between the cluster expansion of \cite{osterwalderseiler1978} and our cluster expansion, is that the base measure off of which we perturb is not simply a product Haar measure, but a different product measure which is not identical edge to edge. The second step is to show that conditional expectations of local observables are almost local which follows from the previously obtained conditional mass gap. Then, finally the third step is based on applying the dynamical approach of \cite{shen2023stochastic} to the model given by the marginal law of the $\mrm{SU}(N)$ field to establish a mass gap for this marginal distribution. Combining all three steps together allows us to recover the mass gap for the $\mrm{U}(N)$ lattice Yang--Mills theory. 

Beyond the immediate problems solved in this paper, we hope that the decomposition of $\mrm{U}(N)$ lattice Yang--Mills as a $\mrm{U}(1) \times \mrm{SU}(N)$ field as in definition \ref{def:random-enviornment}, will continue to be useful in an interplay between results in $\mrm{SU}(N)$ and $\mrm{U}(N)$ lattice Yang--Mills theories.\\

\noindent \textit{Organization:} The rest of the paper is divided into three sections. The bulk of the paper is Section \ref{section:Mass-Gap} where we establish the mass gap result for $\mrm{U}(N)$. This section is broken down into three subsections separating the three key steps of the proof. In subsection \ref{subsection:Step1} we establish a conditional mass gap result for the $\mrm{U}(1)$ field, in subsection \ref{subsection:Step2} we show that conditional observables are very close to being local, and in subsection \ref{subsection:Step3} we show a marginal distribution mass gap for the $\mrm{SU}(N)$ field. The much shorter sections \ref{section:infinite-volume} and \ref{section:Large-N-Limit} establish the infinite volume limit and large N limit results respectively, and heavily rely on the results of section \ref{section:Mass-Gap}.\\

\noindent \textbf{Acknowledgements:} The author would like to thank Sky Cao and Scott Sheffield for many helpful discussions, and for looking over a draft of the manuscript. The author would also like to thank Hao Shen, Rongchan Zhu, and Xiangchan Zhu for helpful email communication, and for pointing out an error in an earlier attempt to prove Theorem \ref{thm:U-N-mass-gap}. R.N. was supported by the NSF under Grant No. GRFP-2141064.

\section{Preliminaries}\label{section:Preliminaries}

\subsection{Lattice Yang--Mills}\label{subsection:Lattice-YM-Prelim}

In this section we recall the lattice Yang--Mills model and set some notation. We will only ever let $G \in \{\mrm{U}(N),\mrm{SU}(N)\}$ for the rest of the paper. The parameter $N \geq 1$ will always refer to the matrix sizes. We will always work with the lattice $\Lambda_L :=[-L,L]^d  \sse \Z^d$. All of our estimates will be uniform in $L$, so we will often just denote $\Lambda_L$ by $\Lambda$.  Let $E_{\Lambda}^+$ (resp. $E_\Lambda$) be the collection of positively oriented (resp. oriented) edges in $\Lambda$. Similarly, let $\mc{P}_{\Lambda}^+$ (resp. $\mc{P}_{\Lambda}$) denote the set of positively oriented (resp. oriented) plaquettes in $\Lambda$. For oriented edges $e \in E_\Lambda$, we denote by $e^{-1}$ the oppositely oriented version of $e$. Similarly, for oriented plaquettes $p \in \mc{P}_\Lambda$, we denote by $p^{-1}$ the oppositely oriented version of $p$.

\begin{definition}[Lattice gauge configuration]\label{def:gauge-configuration}
A lattice gauge configuration $Q$ is a function $Q : E_\Lambda^+ \ra G$. We always implicitly extend $Q : E_\Lambda \ra G$ to oriented edges, by imposing that $Q_{e^{-1}} := Q_e^{-1}$ for all $e \in E_\Lambda^+$.
\end{definition}

\begin{remark}
    We can, and usually will, equivalently identify the set of gauge configurations by $G^{E_{\Lambda}^+}$.
\end{remark}

\begin{definition}[Plaquette variable]
Given a lattice gauge configuration $Q : E_\Lambda^+ \ra \unitary(N)$ and an oriented plaquette $p = e_1 e_2 e_3 e_4$, we define the plaquette variable (abusing notation)
\begin{equs}
Q_p := Q_{e_1} U_{e_2} U_{e_3} U_{e_4}.
\end{equs}
\end{definition}

\begin{remark}\label{rmk:U(1) plaquette}
    In the $\mrm{U}(1)$ case, for each $e$ we can write $Q_e=e^{i\theta_e}$ for some angle $\theta_e$ with $\theta_{e^{-1}}=-\theta_e$, and so $Q_p = e^{i \theta_p}$ where $\theta_p = \theta_{e_1}+\theta_{e_2}+\theta_{e_3}+\theta_{e_4}$.
\end{remark}

\begin{definition}[Orientation sign function]
    For an edge $e$ and a plaquette $p=e_1e_2e_3e_4$, we define the function,
    \begin{equs}
        \mrm{sgn}(e,p):= 
        \begin{cases}
            0 \hspace{4mm} \text{ if } e,e^{-1} \notin \{e_1,e_2,e_3,e_4\},\\
            1 \hspace{4mm} \text{ if } e \in \{e_1,e_2,e_3,e_4\},\\
            -1 \hspace{4mm} \text{ if } e^{-1} \in \{e_1,e_2,e_3,e_4\}.
        \end{cases}
    \end{equs}
\end{definition}
\begin{remark}
    Observe that $\theta_p$ as defined in remark \ref{rmk:U(1) plaquette}, can be written as $\theta_p = \sum_{e \in E_{\Lambda}^+}\mrm{sgn}(e,p)\theta_e$.
\end{remark}

Next, we define the lattice Yang--Mills measure.

\begin{definition}[Lattice Yang--Mills]\label{def:lattice-ym}
The \textit{lattice Yang-Mills measure} with Wilson action, 't Hooft scaling, and gauge group $G \in \{ \mrm{U}(N),\mrm{SU}(N)\}$ is the measure on $G^{E_{\Lambda}^+}$ given by
\begin{equs}
    d\mu_{G,\Lambda, \beta}(Q) := \frac{1}{Z_{G,\Lambda, \beta}} \exp(S_{G,\Lambda,\beta}(Q)) dQ,
\end{equs}
where the action is given by
\begin{equs}
    S_{G,\Lambda,\beta}(Q):=\sum_{p \in \mc{P}_{\Lambda}^+}N\beta \mathrm{Re}\mathrm{Tr}(Q_p)
\end{equs}
and $dQ = \prod_{e \in E_\Lambda^+} dQ_e$, and each $dQ_e$ is an independent copy of Haar measure on $G$ for each edge $e \in E_{\Lambda}^+$. Finally, the partition function,
\begin{equs}
Z_{\mrm{U}(N),\Lambda, \beta} :=\int_{G^{E_{\Lambda}^+}}\prod_{p \in \mc{P}_{\Lambda}^+} \exp(N\beta \mathrm{Re}\mathrm{Tr}(Q_p)) dQ.
\end{equs}
We will sometimes use the notation $\langle \cdot \rangle_{\mu_{G,\Lambda, \beta}}$ or $\langle \cdot \rangle_{G,\Lambda, \beta}$ to denote the expectation with respect to the probability measure $\mu_{G,\Lambda, \beta}$.
\end{definition}

Next, we define a collection of quantities associated to the lattice Yang--Mills measure. These quantities are associated to loops and strings, which we first define.

\begin{definition}[Loops, Loop Variables, and Wilson Loop Observables]\label{def:loops-loop variables}
We will represent loops $\ell$ in $\Lambda$ by the sequence of oriented edges $\ell = e_1 \cdots e_n$ that are traversed by $\ell$. We denote $|\ell| := n$. Let $Q : E_\Lambda^+ \ra G$ be a lattice gauge configuration. For  a loop $\ell=e_1e_2\dots e_n$, let $Q_{\ell}:=Q_{e_1}Q_{e_2}\dots Q_{e_n}$ denote the corresponding loop variable. We define
\begin{equs}
W_\ell(Q) := \tr(Q_\ell),
\end{equs}
where $\tr = \frac{1}{N} \Tr$ is the normalized trace.
We refer to $W_\ell$ as a Wilson loop observable.
\end{definition}



We also set notation regarding local observables, which will be used at various points in the paper.

\begin{notation}[Local observables]\label{notation:local-observables}
A local observable is a function $f \colon G^{\Lambda_f} \ra \R$, where $\Lambda_f \sse E_{\Z^d}^+$ is a finite set. In particular, whenever we specify a local observable $f$, we implicitly also specify the set $\Lambda_f$. Additionally, for any lattice $\tilde{\Lambda}$ such that $\Lambda_f \sse \tilde{\Lambda}$, we may extend $f$ to a function $f \colon G^{E_{\tilde{\Lambda}}^+} \ra \R$.
\end{notation}

\begin{definition}
Given a local observable $f$, we define the following norms:
\begin{equs}
\|f\|_{L^\infty} &:= \sup_{Q \in G^{\Lambda_f}} |f(Q)|, \\
\vertiii{f}_{\infty} &:= \sum_{e \in \Lambda_f }\|\nabla_e f(Q)\|_{L^\infty}.
\end{equs}
Here, $\nabla_e f$ is the gradient of $f$ in the $e$ coordinate, see e.g. \eqref{eq:Edge-grad}.
\end{definition}

\subsection{Riemannian Geometry and Lie Groups}
In this section we set some geometric notation and review the definitions of the geometric objects which will appear in sections \ref{subsection:Step2} and \ref{subsection:Step3}. We will almost always use the same notation as \cite[Section 2]{shen2023stochastic}, nevertheless we review this notation for the sake of completeness. For more details we refer the reader to \cite{RiemannianGeometry}.

\subsubsection{Ricci Curvature and Hessian}
 
For a Riemannian manifold $M$, let $\nabla$ be the Levi-Civita connection associating to vector fields $X$ and $Y$ a vector field $\nabla_Y X$. In our setting, the Levi-Civita will always have a very explicit expression \eqref{eq:Levi-Civita-simplify}.

For $f \in C^\infty(M)$, we denote by $\nabla f$ the gradient vector field of $f$. We also write $\mathrm{Hess}(f)$ for the Hessian. It can be calculated in the following ways:
\begin{equs}\label{eq:Hessian-Formula}
\mathrm{Hess}_f(X,Y) := \mathrm{Hess}(f)(X,Y) = \langle \nabla_X \nabla f, Y \rangle = X(Yf) - (\nabla_X Y)f.
\end{equs}
Since the Levi-Civita connection is torsion-free, $\mathrm{Hess}(f)$ is symmetric in $X,Y$.

The Riemann curvature tensor $\mathscr{R}(\cdot,\cdot)$ associated to vector fields $X,Y$ is an operator defined by
\[
\mathscr{R}(X,Y)Z = \nabla_X (\nabla_Y Z) - \nabla_Y (\nabla_X Z) - \nabla_{[X,Y]} Z.
\]
Let $\{W_i\}_{i=1}^d$ be an orthonormal basis of $T_x M$. The Ricci curvature tensor is defined by
\begin{equs}
\mathrm{Ric}(X,Y) = \sum_{i=1}^d \langle \mathscr{R}(X, W_i) W_i, Y \rangle_{T_x M},
\end{equs}
and is independent of the choice of $\{W_i\}$. Note that $\mathrm{Ric}(X,Y)(x)$ depends on $X,Y$ only via $X(x), Y(x)$ for $x \in M$.

\subsubsection{Lie Groups and Lie Algebras}

For any matrix $A$ we write $A^*$ for the conjugate transpose of $A$. Let $M_N(\mathbb{R})$ and $M_N(\mathbb{C})$ be the space of real and complex $N \times N$ matrices.

For Lie groups $\mathrm{U}(N)$, $\mathrm{SU}(N)$, we write the corresponding Lie algebras as $\mathfrak{u}(N)$, $\mathfrak{su}(N)$ respectively. Every matrix $Q$ in one of these Lie groups satisfies $Q Q^* = I_N$, and every matrix $X$ in one of these Lie algebras satisfies $X + X^* = 0$. Here $I_N$ denotes the identity matrix.

We endow $M_N(\mathbb{C})$ with the Hilbert–Schmidt inner product
\begin{equs}
\langle X, Y \rangle = \mathrm{Re} \, \mathrm{Tr}(X Y^*) \quad \forall X, Y \in M_N(\mathbb{C}).
\end{equs}
We restrict this inner product to our Lie algebra $\mathfrak{g}$. For $X,Y \in \mathfrak{u}(N)$ or $\mathfrak{su}(N)$ we have $\langle X, Y \rangle = -\mathrm{Tr}(XY)$. Note that $\mathrm{Tr}(XY) \in \mathbb{R}$ since $\mathrm{Tr}((XY)^*) = \mathrm{Tr}(Y^* X^*) = \mathrm{Tr}(XY)$, and $\mathrm{Tr}(A^*) = \overline{\mathrm{Tr}(A)}$ for any $A \in M_N(\mathbb{C})$. Throughout the paper,
\begin{equs}
    |X|=\langle X,X\rangle^{1/2}
\end{equs}
will denote the norm with respect to the Hilbert–Schmidt inner product.

In what follows, we always take  $G \in \{\mathrm{U}(N),\mathrm{SU}(N)\}$. Every $X \in \mathfrak{g}$ induces a right-invariant vector field $\widetilde{X}$ on $G$, and for each $Q \in G$, $\widetilde{X}(Q)$ is just given by $XQ$ since $G$ is a matrix Lie group. The inner product on $\mathfrak{g}$ induces an inner product on the tangent space at every $Q \in G$ via the right multiplication on $G$. Hence, for $X,Y \in \mathfrak{g}$, we have $XQ, YQ \in \mrm{T}_Q G$, and their inner product is given by $\mathrm{Tr}((XQ)(YQ)^*) = \mathrm{Tr}(XY^*)$. This yields a bi-invariant Riemannian metric on $G$.

For any function $f \in C^\infty(G)$ and $X \in \mathfrak{g}$, the right-invariant vector field $\widetilde{X}$ induced by $X$ acts on $f$ at $Q \in G$ by the right-invariant derivative
\begin{equs}
\widetilde{X} f(Q) = \left. \frac{d}{dt} \right|_{t=0} f(e^{tX} Q).
\end{equs}
We also have
\begin{equs}\label{eq:commutator-Lie-Group}
[\widetilde{X}, \widetilde{Y}] = \widetilde{[X,Y]}, \quad \text{namely,} \quad ([X,Y] f)(Q) = [XQ, YQ] f(Q),
\end{equs}
where $[\cdot,\cdot]$ is the Lie bracket on $\mathfrak{g}$ on the LHS and the commutator of vector fields on the RHS. Also, for the Levi-Civita connection $\nabla$ we have
\begin{equs}\label{eq:Levi-Civita-simplify}
\nabla_{\widetilde{X}} \widetilde{Y} = \frac12 \widetilde{[X,Y]}.
\end{equs}

\subsubsection{Products of Riemannian Manifolds}
For Riemannian manifolds $M_1, M_2$, the tangent space of the product manifolds $\mrm{T}_{(x_1,x_2)}(M_1 \times M_2)$ is isomorphic to $\mrm{T}_{x_1}M_1 \oplus \mrm{T}_{x_2}M_2$ which is endowed with the inner product
\[
\langle u_1 + u_2, v_1 + v_2 \rangle_{T_{(x_1,x_2)}(M_1 \times M_2)} 
= \langle u_1, v_1 \rangle_{\mrm{T}_{x_1} M_1} + \langle u_2, v_2 \rangle_{\mrm{T}_{x_2} M_2}.
\]
For a finite collection of Riemannian manifolds $(M_u)_{u \in A}$ where $A$ is some finite set, the product is defined analogously.

If all $M_e$ are the same manifold $M$, the product is written as $M^A$. In this case, given a point $x = (x_e)_{e \in A} \in M^A$, if $v_e \in \mrm{T}_{x_e} M_e$ for some $x_e \in M_e$, we will sometimes view $v_e$ as a tangent vector in $\mrm{T}_x M^A$ which has zero components for all $e' \neq e$. Continuing with this notation, if $\{ v_e^i \}_{i=1}^d$ is a basis (resp. orthonormal basis) of $T_{x_e} M_e$, then $\{ v_e^i \}_{e \in A, \, i=1,\dots,d}$ is a basis (resp. orthonormal basis) of $\mrm{T}_x M^A$.

For Lie groups $G_1, G_2$, the group multiplication is defined on $G_1 \times G_2$ componentwise. The Lie algebra $\mathfrak{g}$ of $G_1 \times G_2$ is isomorphic to $\mathfrak{g}_1 \oplus \mathfrak{g}_2$ where $\mathfrak{g}_i$ is the Lie algebra of $G_i$. The Lie bracket on $\mathfrak{g}_1 \oplus \mathfrak{g}_2$ is defined componentwise. If $X = (X_1, X_2) \in \mathfrak{g}$, then the induced right-invariant vector field $\widetilde{X}(x)$ for every $x \in G_1 \times G_2$ is equal to $(\widetilde{X}_1(x), \widetilde{X}_2(x))$. In particular, \eqref{eq:commutator-Lie-Group} and \eqref{eq:Levi-Civita-simplify} still hold for any two right-invariant vector fields on the Lie group product.

With similar notation as above we can define the product $G^A$ and its Lie algebra $\mathfrak{g}^A$ for a finite set $A$. Given $X \in \mathfrak{g}^A$, the exponential map $t \mapsto \exp(tX)$ is also defined pointwise as
\[
\exp(tX)_e := e^{t X_e} \quad \text{for each } e \in A.
\]
In the following we consider the configuration space $G^{E_{\Lambda_L}^+}$ which is the product Lie group with $G \in \{\mrm{U}(N),\mrm{SU}(N)\}$
consisting of all maps $Q : e \in E_{\Lambda_L}^+ \mapsto Q_e \in G$. The corresponding Lie algebra is $\mathfrak{g}^{E_{\Lambda_L}^+}$, the direct sum of $\mathfrak{g}$. For any matrix-valued functions $A,B$ on $V_{\Lambda_L}$, we denote by $AB$ the pointwise product $(A_e B_e)_{u \in E_{\Lambda_L}^+}$.

As above, the tangent space at $Q \in G^{E_{\Lambda}^+}$ consists of the products $XQ = (X_u Q_u)_{u \in V_{\Lambda_L}}$ with $X \in \mathfrak{g}^{E_{\Lambda_L}^+}$, and given two such elements $XQ$ and $YQ$, their inner product is defined by
\[
\langle XQ, YQ \rangle_{\mrm{T}_Q G^{E_{\Lambda_L}}} = \sum_{u \in V_{\Lambda_L}} \mathrm{Tr}(X_u Y_u^*).
\]
A basis of the tangent space $\mrm{T}_Q G^{E_{\Lambda_L}}$ is given by $\{ X_e^i Q : e \in E_{\Lambda_L^+}, \ 1 \le i \le d(\mathfrak{g}) \}$ where for each $e$, $\{ X_e^i \}_i$ is a basis for $\mathfrak{g}$.

Given any function $f \in C^\infty(G^{E_{\Lambda_L}^+})$, the right-invariant derivative is given by
\[
\widetilde{X} f(Q) = \left. \frac{d}{dt} \right|_{t=0} f(\exp(tX) Q).
\]
For each $Q \in G^{E_{\Lambda_L}^+}$, the gradient $\nabla f(Q)$ is the element of the tangent space at $Q$ which satisfies, for each $X \in \mathfrak{g}^{E_{\Lambda_L}^+}$,
\begin{equation}
\langle \nabla f(Q), XQ \rangle_{\mrm{T}_Q G^{E_{\Lambda_L}^+}} = \left. \frac{d}{dt} \right|_{t=0} f(\exp(tX) Q) = (XQ) f.
\end{equation}
We can write
\[
\nabla f = \sum_{i=1}^{d(\mathfrak{g})} \sum_{u \in E_{\Lambda_L}^+} (v_e^i f) \, v_e^i
\]
with $\{ v_e^i : e \in E_{\Lambda_L}^+, \ i=1,\dots,d(\mathfrak{g}) \}$ being an orthonormal basis of $\mrm{T}_Q G^{E_{\Lambda_L}^+}$. We then define
\begin{equs}\label{eq:Edge-grad}
    \nabla_e f \overset{\mathrm{def}}{=} \sum_{i=1}^{d(\mathfrak{g})} (v_e^i f) \, v_e^i,
    \qquad
    \Delta_e f := \mathrm{div} \, \nabla_e f = \sum_{i=1}^{d(\mathfrak{g})} \langle \nabla_{v_e^i} \nabla_e f, v_e^i \rangle.
\end{equs}

Here $\nabla_e$ and $\Delta_e$ can be viewed as the gradient and the Laplace–Beltrami operator (with respect to the variable $Q_e$) on $G$ endowed with the metric given above.

\subsubsection{Lie Group Brownian Motion}
Denote by $\mf{B}$ and $B$ the Brownian motions on a Lie group $G$ and its Lie algebra $\mathfrak{g}$ respectively. The Brownian motion $B$ is characterized by
\begin{equs}\label{eq:Lie-Algebra-BM-cov}
\mathbb{E} \left[ \langle B(s), X \rangle \langle B(t), Y \rangle \right] = \min(s,t) \langle X, Y \rangle, \quad \forall X, Y \in \mathfrak{g}.
\end{equs}
By \cite[Sec.~1.4]{Levy2011a}, the Brownian motions $\mf{B}$ and $B$ are related through the following SDE:
\begin{equation}
d\mf{B} = dB \circ \mf{B} = dB \, \mf{B} + \frac{c_{\mathfrak{g}}}{2} \mf{B} \, dt,
\end{equation}
where $\circ$ is the Stratonovich product, and $dB \, \mf{B}$ is in the Itô sense. Here the constant $c_{\mathfrak{g}}$ is determined by $\sum_\alpha v_\alpha^2 = c_{\mathfrak{g}} I_N$, where $\{ v_\alpha \}_{\alpha=1}^{d(\mathfrak{g})}$ is an orthonormal basis of $\mathfrak{g}$. Moreover, by \cite[Lem.~1.2]{Levy2011a},
\begin{equation}
c_{\mathfrak{su}(N)} = -\frac{N^2 - 1}{N}.
\end{equation}

\subsection{$\mrm{U}(N)$ Yang-Mills as random environment $\mrm{SU}(N)$ Field}

We will always identify elements $z \in \mrm{U}(1)$ with $z = e^{i \theta}$ for an angle $\theta \in [0,2\pi)$.

\begin{definition}\label{def:random-enviornment}
    Take $\mc{U}_L=\mc{U}$ to be the space $[0,2\pi)^{\Lambda_L} \times \mrm{SU}(N)^{\Lambda_L}$, and define the measure $\mu_{\mc{U}_L,\beta}$ on $\mc{U}_L$ by setting,
    \begin{equs}
        S_{\mc{U}_L}(\theta,Q)=S_{\mc{U}}(\theta,Q):=N\beta \sum_{p \in \mc{P}_{\Lambda_L}^+}\mrm{Re}(e^{i \frac{\theta_p}{N}}\mrm{Tr}(Q_p))
    \end{equs}
    \begin{equs}
        d\mu_{\mc{U}_L,\beta}(\theta,Q)=d\mu_{\Lambda_L,\beta}(\theta,Q) := \frac{1}{Z_{\mc{U}_L,\beta}} \exp(S_{\mc{U}}(\theta,Q)) d\theta dQ
    \end{equs}
    Where $d\theta = \prod_{e \in E_{\Lambda}^+} \mathbbm{1}_{0 \leq \theta < 2\pi}d\theta_e$ and $dQ = \prod_{e \in E_{\Lambda}^+} dQ_e$ where each $dQ_e$ is an $\mrm{SU}(N)$ Haar measure.
\end{definition}

\begin{lemma}\label{lm:measure-decompostion}
    Suppose the field $(\theta,Q)=(\theta_e,Q_e)_{e \in \Lambda}$ is a random field with Law given by $\mu_{\mc{U}_L,\beta}$. Then $(e^{i \frac{\theta_e}{N}}Q_e)_{e \in E_{\Lambda}^+}$ is distributed according to the $\mrm{U}(N)$ lattice Yang-Mills measure on $\Lambda_L$ with coupling constant $\beta$, $\mu_{\mrm{U}(N),\Lambda_L,\beta}$.
\end{lemma}

\begin{proof}
    The proof follows from the fact that if $\theta$ is uniformly distributed on $[0,2\pi)$, and $Q$ is independent of $\theta$, and distributed according to the Haar measure on $\mrm{SU}(N)$, then $e^{i \frac{\theta}{N}}Q$ is distributed according to the Haar measure on $\mrm{U}(N)$ (i.e. see \cite[Lemma 5.15]{CompactLieGroupRMT}). 
\end{proof}

\section{Mass Gap}\label{section:Mass-Gap}
Throughout this section we use $C,C_0,C_1,C_2,\dots$ to denote unspecified constants, and when we write $C=C(a,b,c,\dots)$, we mean the constant $C$ depends on the parameters $a,b,c,\dots$. Sometimes line to line, constants with the same name will have different values. Lastly we emphasize that none of the constants in this section depend on $L$.

In order to prove the mass gap statement of equation \eqref{eq:main-mass-gap}, we will need to decompose the $\mrm{U}(N)$ covariance into three parts which will each be bounded separately in the next three subsections. The decomposition of the covariance is based on the $\mrm{U}(N) \simeq \mrm{U}(1) \times \mrm{SU}(N)$ measure decomposition of Lemma \ref{lm:measure-decompostion}. We now state the three propositions corresponding to these three parts. In the following propositions we set $\beta^*=\beta^*(d):=10^{-6d}$ and $\tilde{\beta}=\tilde{\beta}(d)$ will be some fixed constant possibly smaller than $\beta^*(d)$. Additionally, all expectations in the propositions and remainder of Section \ref{section:Mass-Gap} are assumed to be taken with respect to $\mu_{\mc{U}_L,\beta}$ unless otherwise specified.

The statement of the first of the three propositions below is stating a mass gap for the conditional law of $\mu_{\mc{U}_L,\beta}$ on $\mrm{U}(1)^{E_{\Lambda}^+}$ which is obtained by fixing the $\mrm{SU}(N)^{E_{\Lambda}^+}$ part of the field (i.e. conditioning $Q=Q' \in \mrm{SU}(N)^{E_{\Lambda}^+}$).

\begin{prop}\label{prop:conditional-mass-gap}
    Suppose $N>8\pi$ and $\beta<\beta^*(d)=10^{-6d}$ . Fix an $\mrm{SU}(N)$ field $(Q_e')_{e \in E_{\Lambda}^+}$. Then for any two local observables $f,g \in C^{\infty}(\mc{U}_L)$,
    \begin{equs}
        |\E[fg| Q = Q']-\E[f| Q = Q']\E[g| Q = Q']| \leq C_1 \|f\|_{L^{\infty}}\|g\|_{L^{\infty}}10^{-2d \cdot d(\Lambda_f,\Lambda_g)}
    \end{equs}
     for a universal constant  $C=C(d,|\Lambda_f|,|\Lambda_g|)$.
\end{prop}

The next proposition roughly states that when a local observable $f \in C^{\infty}(\mc{U}_L)$ is averaged over the $\theta$ field, holding the $Q$ field fixed, the conditional expectation which is now a function of only $Q \in \mrm{SU}(N)^{E_{\Lambda}^+}$ is still almost a local observable.

\begin{prop}\label{prop:Localizing conditional observables}
    Suppose $\beta<\beta^*(d)=10^{-6d}$, $N>8\pi$ and fix any local observable $f \in C^{\infty}(\mc{U}_L)$ and an $\mrm{SU}(N)$ field $(Q_e')_{e \in E_{\Lambda}^+}$. Then for any fixed distance $R>0$, there are constants $C_1=C_1(d,N,|\Lambda_f|)$ and $C_2=C_2(d)$  such that,
    \begin{equs}\label{eq:one-func-sensativity}
        |\E[f|Q=Q']-\E[f|Q_e=Q'_e \text{ for } d(e,\Lambda_f) \leq R]]| \leq C_1 \|f\|_{L^{\infty}}e^{-C_2 R},
    \end{equs}
    and as a corollary, for two local observables $f,g \in C^{\infty}(\mc{U}_L)$ there are constants $C_1'=C_1'(d,N,|\Lambda_f|,|\Lambda_g|)$ and $C_2'=C_2'(d)$ such that,
    \begin{equs}\label{eq:two-func-sensativity}
        ~&|\E[f|Q]\E[g|Q]-\E[f|(Q_e)_{e \in \Lambda_{f,g^c}}]\E[g|(Q_e)_{e \in \Lambda_{f^c,g}}]| \\&\leq C_1' \|f\|_{L^{\infty}} \|g\|_{L^{\infty}} e^{-C_2' d(\Lambda_f,\Lambda_g)},
    \end{equs}
    where $\Lambda_{f,g^c}:=\{e:d(e,\Lambda_f) \leq \frac{1}{3}d(\Lambda_f,\Lambda_g)\}$ and $\Lambda_{f^c,g}:=\{e:d(e,\Lambda_g) \leq \frac{1}{3}d(\Lambda_f,\Lambda_g)\}$.
\end{prop}

In order to state the third proposition, we introduce the measure $\tilde{\mu}_{\mrm{SU}(N),\Lambda_L,\beta}$ which is the marginal distribution of the $\mrm{SU}(N)$ part of $\mu_{\mc{U}_L,\beta}$. It will be important for us that we can express $\tilde{\mu}_{\mrm{SU}(N),\Lambda_L,\beta}$ as a measure on $\mrm{SU}(N)^{\Lambda_L}$ given by,
\begin{equs}\label{eq:SU(N)-marginal-dist}
    d\tilde{\mu}_{\mrm{SU}(N),\Lambda_L,\beta}(Q):=\frac{1}{\tilde{Z}_{\mrm{SU}(N),\Lambda_L,\beta}} \exp(\tilde{S}(Q)) dQ, 
\end{equs}
where $dQ = \prod_{e \in E_{\Lambda}^+}$ is the product Haar measure and
\begin{equs}\label{eq:marginal-action}
    \tilde{S}(Q):=\log\bigg( \int \exp(S_{\mc{U}}(\theta,Q))d\theta\bigg),
\end{equs}
and we recall that $S_{\mc{U}}(\theta,Q):= N\beta \sum_{p \in \mc{P}_{\Lambda}^+} \mrm{Re}(e^{i\frac{\theta_p}{N}} \mrm{Tr}(Q_p))$.

The final proposition is a mass gap for this marginal law.
\begin{prop}\label{prop:averaged-out-mass-gap}
    There exists a constant $\tilde{\beta}=\tilde{\beta}(d)$ such that whenever $\beta<\tilde{\beta}$, $N>8\pi$ , there exist constants $C_1,C_2$ only depending on $N$, and $d$, such that for all local observables $f,g\in C^{\infty}(\mrm{SU}(N)^{E_{\Lambda}^+})$,
    \begin{equs}
        \mrm{Cov}_{\tilde{\mu}_{\mrm{SU}(N),\Lambda_L,\beta}}(f,g) \leq C_1(\|f\|_{L^2}\|g\|_{L^2}+\vertiii{f}_{\infty} \vertiii{g}_{\infty})e^{-C_2 d(\Lambda_f, \Lambda_g)}.
    \end{equs}
\end{prop}
The desired mass gap \eqref{eq:main-mass-gap} readily follows from these three propositions.

\begin{proof}[Proof of Mass Gap \eqref{eq:main-mass-gap}]
    From the law of total expectation, we can rewrite,
    \begin{equs}
        \mrm{Cov}_{\mu_{\mrm{U}(N),\Lambda_L,\beta}}(f,g) &= \E[fg]-\E[f]\E[g]\\
        &= (\E[\E[fg|Q]-\E[f|Q]\E[g|Q]])\\
        &+(\E[\E[f|Q]\E[g|Q]-\E[f|(Q_e)_{e \in\Lambda_{f,g^c}}]\E[g|(Q_e)_{e \in\Lambda_{f^c,g}}]])\\
        &+(\E[\E[f|(Q_e)_{e \in\Lambda_{f,g^c}}]\E[g|(Q_e)_{e \in\Lambda_{f^c,g}}]]-\E[\E[f|(Q_e)_{e \in\Lambda_{f,g^c}}]]\E[\E[g|(Q_e)_{e:e \in\Lambda_{f^c,g}}]]),
    \end{equs}
    where expectations/conditional expectations without reference to a probability measure are assumed to be taken with respect to $\mu_{\mc{U}_L,\beta}$, and as in Proposition \ref{prop:Localizing conditional observables}, $\Lambda_{f,g^c}:=\{e:d(e,\Lambda_f) \leq \frac{1}{3}d(\Lambda_f,\Lambda_g)\}$ and $\Lambda_{f^c,g}:=\{e:d(e,\Lambda_g) \leq \frac{1}{3}d(\Lambda_f,\Lambda_g)\}$.
    Now each of the three terms in the decomposition above can be bounded by $C_1 \big(\|f\|_{L^{\infty}}+\vertiii{f}_{\infty}\big)\big(\|g\|_{L^{\infty}}+\vertiii{g}_{\infty}\big) e^{-C_2 d(\Lambda_f,\Lambda_g)}$, by applying Propositions \ref{prop:conditional-mass-gap}, \ref{prop:Localizing conditional observables}, and \ref{prop:averaged-out-mass-gap} respectively. The proofs of these propositions is the content of the next three subsections. The $N>8\pi$ is unimportant since the result of \cite{osterwalderseiler1978} applies when $N \leq 8\pi$ for sufficiently small $\beta$.
\end{proof}

Before providing the full details in the next three subsections, we give a sketch of the proofs below.\\

\textit{Sketch of Step 1 (Proposition \ref{prop:conditional-mass-gap}):} The key lemma here is Lemma \ref{lm:Cluster-Exp-Decomp} which will allow us to write the conditional law $\mu_{\mc{U}}(\theta|Q=Q')$ in a form which is suitable to perform a cluster expansion. The idea behind Lemma \ref{lm:Cluster-Exp-Decomp} is that if we fix $Q_e=Q'_e$ for all $e$, and treat the action $S_{\mc{U}}(\theta,Q')$ as a function of only $\theta$, then replacing the action $S_{\mc{U}}(\theta,Q')$ by its linear approximation (1st order Taylor series), leads to a new measure $\nu_{Q'}(\theta)$ which separates as a product of measures on each edge variable. We will treat $\mu_{\mc{U}}(\theta|Q=Q')$ as a perturbation of $\nu_{Q'}(\theta)$, performing a cluster expansion identical to that in \cite{osterwalderseiler1978}.\\

\textit{Sketch of Step 2 (Proposition \ref{prop:Localizing conditional observables}):} This step is the easiest. If we simply differentiate the function $\phi(Q')=\E[f|Q=Q']$ with respect to an edge variable $Q_e'$ for an edge $e$ far away from $\Lambda_f$, then the computation will lead us to an answer which is a conditional covariance between $f$ an observable around $e$. Thus plugging in the conditional mass gap obtained in step 1, gives us an exponential bound in $d(e,\Lambda_f)$ for this derivative. Proposition \ref{prop:Localizing conditional observables} then essentially follows from the fundamental theorem of calculus.\\

\textit{Sketch of Step 3 (Proposition \ref{prop:averaged-out-mass-gap}):} The final step is based on the dynamical approach of \cite{shen2023stochastic}. We now use the action $\tilde{S}(Q)$ for the marginal law $\tilde{\mu}_{\mc{U}_L,\beta}$ to define the Langevin dynamics on $\mrm{SU}(N)^{E_{\Lambda}^+}$ which we informally write as
\begin{equs}
    dQ_t = \nabla\tilde{S}(Q_t) dt+\sqrt{2}dB_t.
\end{equs}
The first step is to use the Bakry-Emery condition to verify that this dynamics mixes exponentially fast to $\tilde{\mu}_{\mc{U}_L,\beta}$. The necessary computation involves bounding $\mrm{Hess}(\tilde{S})(v,v)$ which turns out to rely on the corresponding argument in \cite{shen2023stochastic} together with the conditional mass gap of step 1. The next step is to verify that the Langevin dynamic we study in this step is ``local'' in a sense which is made precise in Lemma \ref{lm:Generator-Commutator-Bound}. Piecing together the mixing and locality properties, the rest of the argument is a modification of the arguments presented in \cite[Section 8.3]{GuionnetLogSobolev} with an additional combinatorial bound.

\subsection{Step 1: Mass gap conditioned on $\mrm{SU}(N)$ field}\label{subsection:Step1}

The main result of this section is the following proposition.

In this section we set $\beta^*(d)=10^{-6d}$, and always assume $\beta \in (0,\beta^*(d))$.

The key idea in this section is the following representation of $d\mu(\theta,Q)$ and $d\mu(\theta|Q=Q')$ which allows for a cluster expansion.

\begin{lemma}\label{lm:Cluster-Exp-Decomp}
    We can express $d\mu_{\mc{U}_L,\beta}(\theta,Q)$ as
    \begin{equs}\label{eq:dmu-cluster-exp}
        d\mu_{\mc{U}_L,\beta}(\theta,Q) = \frac{1}{Z_{\mc{U}_L,\beta}} \prod_{p \in \mc{P}_{\Lambda}^+}(1+\varphi_p) \prod_{e \in E_{\Lambda}^+} e^{- \beta \sum_{p \in \mc{P}_{\Lambda}^+} \mrm{sgn}(e,p) \mrm{Im}\mrm{Tr}(Q_p) \theta_e} d\theta_e d\mu_{\mrm{SU}(N),\Lambda_L,\beta},
    \end{equs}
    where $\varphi_p := \varphi(\theta_p,Q_p)$ for $\varphi(\theta,Q) := e^{N\beta \mrm{Re}((e^{i\frac{\theta}{N}}-1-i\frac{\theta}{N})\mrm{Tr}(Q))}-1$, and for any edge $e$ and plaquette $p$, $\mrm{sgn}(e,p)=1$ if $e$ appears in $p$, $\mrm{sgn}(e,p)=-1$ if $e^{-1}$ appears in $p$, and $\mrm{sgn}(e,p)=0$ otherwise. 

    As a consequence we can express the conditional distribution $d\mu_{Q'}(\theta):=d\mu_{\mc{U}_L,\beta}(\theta,Q | Q=Q')$ as,
    \begin{equs}
        d\mu_{Q'}(\theta):=\frac{1}{Z(Q')} \prod_{p \in \mc{P}_{\Lambda}^+}(1+\varphi(\theta_p,Q_p')) \prod_{e \in E_{\Lambda}^+} d\nu_e(\theta_e). \label{eq:conditional-dmu-cluster-exp}
    \end{equs}
    where
    \begin{equs}
        d\nu_e(\theta_e):=\frac{1}{Z_e(Q')}e^{- \beta \sum_{p \in \mc{P}_{\Lambda}^+} \mrm{sgn}(e,p) \mrm{Im}\mrm{Tr}(Q_p) \theta_e} d\theta_e,
    \end{equs}
    for 
    \begin{equs}
        Z_e(Q'):=\int_{0}^{2\pi}e^{- \beta \sum_{p \in \mc{P}_{\Lambda}^+} \mrm{sgn}(e,p) \mrm{Im}\mrm{Tr}(Q_p) \theta_e} d\theta_e,
    \end{equs}
     and where 
     \begin{equs}
         Z(Q'):=\int \prod_{p \in \mc{P}_{\Lambda}^+}(1+\varphi(\theta_p,Q_p')) \prod_{e \in E_{\Lambda}^+} d\nu_e(\theta_e).
     \end{equs}

    Finally observe that if $N >8\pi$ and $\beta \leq \beta^*(d)$, then $|\varphi(\theta_p,Q_p)| \leq  10^{4-6d}$ for any $p$ and choice of fields $(\theta,Q)$.
\end{lemma}

\begin{proof}
    In order to prove \eqref{eq:dmu-cluster-exp}, note that
    \begin{equs}
        S_{\mc{U}}(\theta,Q) &= N\beta \sum_{p \in \mc{P}_{\Lambda}^+} \mrm{Re}\,e^{i \frac{\theta_p}{N}} \mrm{Tr}\,Q_p\\
        &=N\beta \sum_{p \in \mc{P}_{\Lambda}^+} \mrm{Re}\, \left(\left(e^{i \frac{\theta_p}{N}}-1-i\frac{\theta_p}{N}\right)\mrm{Tr}\,Q_p\right)+N\beta \sum_{p \in \mc{P}_{\Lambda}^+} \mrm{Re}\,\left(i\frac{\theta_p}{N}\mrm{Tr}\,Q_p\right)+ +N\beta\sum_{p \in \mc{P}_{\Lambda}^+} \mrm{Re}\, \mrm{Tr}\,Q_p\\
        & = N\beta \sum_{p \in \mc{P}_{\Lambda}^+} \mrm{Re}\, \left(\left(e^{i \frac{\theta_p}{N}}-1-i\frac{\theta_p}{N}\right)\mrm{Tr}\,Q_p\right)-\beta\sum_{e \in E_{\Lambda}^+} \sum_{p \in \mc{P}_{\Lambda}^+}\mrm{sgn}(e,p)\mrm{Im}\,\mrm{Tr}\,Q_p \theta_e + S_{\mrm{SU}(N),\Lambda_L,\beta}(Q),
    \end{equs}
    where the last equality used the relation  $\theta_p = \sum_{e } \mrm{sgn}(e,p)\theta_e$ followed by an interchange of summation.

    As a result,
    \begin{equs}
        e^{S_{\mc{U}}(\theta,Q)}&=\bigg( \prod_{p \in \mc{P}_{\Lambda}^+}e^{N\beta  \mrm{Re} \,((e^{i \frac{\theta_p}{N}}-1-i\frac{\theta_p}{N})\mrm{Tr}\,Q_p)}\bigg) \bigg( \prod_{e \in E_{\Lambda}^+}e^{-\beta \sum_{p \in \mc{P}_{\Lambda}^+}\mrm{sgn}(e,p)\mrm{Im}\,\mrm{Tr}\,Q_p \theta_e}\bigg) e^{S_{\mrm{SU}(N),\Lambda_L,\beta}(Q)}\\
        &=\bigg( \prod_{p \in \mc{P}_{\Lambda}^+}(1+\varphi_p)\bigg) \bigg( \prod_{e \in E_{\Lambda}^+}e^{-\beta \sum_{p \in \mc{P}_{\Lambda}^+}\mrm{sgn}(e,p)\mrm{Im}\,\mrm{Tr}\,Q_p \theta_e}\bigg) e^{S_{\mrm{SU}(N),\Lambda_L,\beta}(Q)},
    \end{equs}
    from which \eqref{eq:dmu-cluster-exp} readily follows. The conditional distribution formula \eqref{eq:conditional-dmu-cluster-exp} is an immediate corollary.
    
    For the bound on $\varphi$ observe that $|\frac{\theta_p}{N}| \leq \frac{8\pi}{N}$, so
    \begin{equs}
        e^{i \frac{\theta_p}{N}} = 1+i \frac{\theta_p}{N}+(e^{i \frac{\theta_p}{N}}-1-i \frac{\theta_p}{N}),
    \end{equs}
    with  
    \begin{equs}
        |e^{i \frac{\theta_p}{N}}-1-i \frac{\theta_p}{N}| &\leq \sum_{k=2}^{\infty} \frac{(8\pi/N)^k}{k!}\\
        &\leq (8\pi/N)^2 \sum_{k=0}^{\infty} (4 \pi /N)^{k}\\
        & \leq \frac{128 \pi^2}{N^2},
    \end{equs}
    where we used the fact that $\frac{4\pi}{N}\leq \frac{1}{2}$. As a consequence if $\beta \leq 10^{-6d}$, then,
    \begin{equs}
        |N\beta \mrm{Re}\,((e^{i\frac{\theta}{N}}-1-i\frac{\theta}{N})\mrm{Tr}(Q))| \leq 128\pi^2 \beta \leq 5000 \beta =  5\cdot 10^{3-6d}.
    \end{equs}
    Thus using the inequality $|e^x-1| \leq 2x$ for $|x|\leq \frac{1}{2}$, we get the desired bound on $\varphi_p$.
\end{proof}

The rest of the section is dedicated to carrying out a cluster expansion argument similar to \cite{osterwalderseiler1978}. We include some of the technical details which are left out of \cite{osterwalderseiler1978}. To carry out the cluster expansion we first need to define some terminology relating to clusters.

\begin{definition}[Graph structure on plaquettes and clusters]
We say that two plaquettes $p,p' \in \mc{P}_{\Lambda}^+$ are neighbors if they share a common edge $e \in E_{\Lambda}$, and we write $p \sim p'$. Given a plaquette set $K \sse \mc{P}_\Lambda^+$ and a subgraph $\Lambda' \subset \Lambda$, we will say that $K$ is a \textit{cluster} of $\Lambda'$ if every connected component of $K$ contains at least one edge of $\Lambda'$. 
\end{definition}

For a subgraph $\Lambda' \subset \Lambda$, we will let $\mc{C}(\Lambda')$ denote the set of clusters of $\Lambda'$, and we let $\mc{C}(f)$ denote the set of connected clusters of $\Lambda_f$. 

For a subgraph $\Lambda' \subset \Lambda$, we write $\Lambda \backslash \overline{\Lambda'}$ to denote the graph with vertices which do not belong to any edges in $\Lambda'$, and all the edges from $E_{\Lambda}$ connecting such vertices. Moreover for any subgraph $\Lambda'$ we set $\mc{P}_{\Lambda'}^+$ to be the set of positively oriented plaquettes in $\Lambda$ for which all of its edges are contained in the edge set of $\Lambda'$, $E_{\Lambda'}$. Lastly we can define the partition function restricted to $\Lambda'$
\begin{equs}
    Z_{\Lambda'}(Q'):=\int \prod_{p \in \mc{P}_{\Lambda'}^+}(1+\varphi(\theta_p,Q'_p)) \prod_{e \in E_{\Lambda'}^+} d\nu_e(\theta_e).
\end{equs}

The cluster expansion is the following representation for conditional expectations which follows from simply expanding the product $\prod_{p \in \mc{P}_{\Lambda'}^+}(1+\varphi(\theta_p,Q'_p))$.

\begin{lemma}
    For any local observable $f$, we can write 
    \begin{equs}
        \E[f(\theta,Q)|Q=Q'] = \sum_{K \in \mc{C}(f)} \int  f(\theta,Q')\prod_{p \in K}\varphi(\theta_p,Q_p')\prod_{e \in E_K} d\nu_e(\theta_e)\frac{Z_{\Lambda \backslash \overline{K \cup \Lambda_f}}(Q')}{Z_{\Lambda}(Q')}. 
    \end{equs}
\end{lemma}

\begin{remark}
    Here we write $E_K^+$ to denote the set of all positively oriented edges which belong to a plaquette in $K$. Additionally, $\overline{K \cup \Lambda_f}$ denotes the subgraph of $\Lambda$ containing all edges which do not belong to a plaquette which is connected to $K \cup \mc{P}_{\Lambda_f}^+$. 
\end{remark}

\begin{remark}
    To avoid cumbersome notation, we write some of the lemmas below, only in the case without boundary conditions, however, these results all extend to the case with boundary.
\end{remark}

The main consequence of the cluster expansion which we need is the following lemma on sensitivity to boundary conditions. In the statement we use $(\Lambda')^c$ to denote $\Lambda \backslash \Lambda'$ and the distance between two subgraphs $d(\cdot,\cdot)$ will denote the usual graph distance.

\begin{lemma}[Sensitivity to boundary conditions]\label{lm-Sensitivity-to-boundary-conditions}
    Suppose $\beta<\beta^*(d)$ and let $f$ be a local observable, and $\Lambda'\subset\Lambda$ be a subgraph containing $\Lambda_f$ with $d(\Lambda_f,(\Lambda')^c)>R$. Then letting for any two fixed boundary conditions $(\theta'_e)_{e \in E_{(\Lambda')^c}^+},(\theta''_e)_{e \in E_{(\Lambda')^c}^+}$, we have,
    \begin{equs}
        |\E[f|Q=Q',(\theta_e)_{e \in E_{(\Lambda')^c}^+}=(\theta'_e)_{e \in E_{(\Lambda')^c}^+}]-\E[f|Q=Q',(\theta_e)_{e \in E_{(\Lambda')^c}^+}=(\theta''_e)_{e \in E_{(\Lambda')^c}^+}]| \leq C_1 \|f\|_{L^{\infty}} 10^{-2d R},
    \end{equs}
    for some universal constant $C_1=C_1(d,|\Lambda_f|)$.
\end{lemma}

Assuming this key consequence of the cluster expansion, we can readily prove the conditional mass gap, Proposition \ref{prop:conditional-mass-gap}.

\begin{proof}[Proof of Proposition \ref{prop:conditional-mass-gap}]
    Let $\Lambda' :=\{v \in \Lambda: d(v,\Lambda_f \cup \Lambda_g) \leq \frac{1}{3} d(\Lambda_f,\Lambda_g)\}$. Then conditioning on the values of $\theta$ in the interface $(\Lambda')^c$, and using the law of total expectation,
    \begin{equs}
        ~&|\E[fg| Q ]-\E[f| Q]\E[g| Q]| \\
        &=|\E[\E[fg|(\theta_e)_{e \in E_{(\Lambda')^c}},Q]|Q]-\E[f|Q]\E[g|Q]|
    \end{equs}
    
    Next since $f$ and $g$ are independent after conditioning on both $(\theta_e)_{e \in E_{(\Lambda')^c}}$ and $Q$, and $\E[g|Q]$ is measurable with respect to the sigma algebra generated by the value of $Q$,
    \begin{equs}
        ~&|\E[\E[fg|(\theta_e)_{e \in E_{(\Lambda')^c}},Q]|Q]-\E[f|Q]\E[g|Q]|\\
        &=|\E[\E[f|(\theta_e)_{e \in E_{(\Lambda')^c}},Q](\E[g|(\theta_e)_{e \in E_{(\Lambda')^c}},Q]-\E[g|Q])Q])|Q]| \\
        &\leq C_1(d,|\Lambda_g|)\|f\|_{\infty}\|g\|_{\infty} 10^{-2d \cdot d(\Lambda_f,\Lambda_g)}.,
    \end{equs}
    where in the last line we applied the conditional Hölder inequality along with Lemma \ref{lm-Sensitivity-to-boundary-conditions}.
\end{proof}

Now we move towards proving Lemma \ref{lm-Sensitivity-to-boundary-conditions}. The first step towards proving Lemma \ref{lm-Sensitivity-to-boundary-conditions} is the following bound on terms corresponding to large clusters.

\begin{lemma}\label{lm:cluster-expansion}
    Suppose $\beta<\beta^*(d)$. For any local observable $f$, and a fixed $R>0$, there is some constant $C=C(d,|\Lambda_f|)$ such that, 
    \begin{equs}\label{eq:Cluster-Expansion-Bound}
        \sum_{K \in \mc{C}(f):|K|>R}\bigg|\int  f(\theta,Q')\prod_{p \in K}\varphi(\theta_p,Q_p')\prod_{e \in E_K} d\nu_e(\theta_e) \frac{Z_{\Lambda \backslash \overline{K \cup \Lambda_f}}(Q')}{Z_{\Lambda}(Q')}\bigg| \leq C 10^{-3d R} \|f\|_{L^{\infty}} 
    \end{equs}
\end{lemma}

The proof of Lemma \ref{lm:cluster-expansion} follows from an estimate on ratios of partition functions and a combinatorial lemma bounding the number of clusters of a fixed size.

\begin{lemma}\label{lm:Comparison-of-partition-functions}
    Suppose $\beta<\beta^*(d)$. Then 
    \begin{equs}
       \frac{Z_{\Lambda \backslash p }(Q')}{Z_{\Lambda}(Q')} \leq (1-10^{4-6d})^{-1}
    \end{equs}
\end{lemma}

\begin{proof}
    Since we are only interested in real $\beta>0$, the proof follows from the trivial bound $\sup_{\theta,Q'} |\varphi(\theta_p,Q'_p)| \leq 10^{4-6d}$.
\end{proof}

\begin{lemma}\label{lm:Combo-Lemma}
    \begin{equs}
        |\{K \in \mc{C}(f): |K|=m\}|\leq e^{2d|\Lambda_f|}(40)^{md}
    \end{equs}
\end{lemma}

\begin{proof}
    The proof is a standard argument which is given in \cite[Lemma 4.6]{CNS2025a}.
\end{proof}

\begin{proof}[Proof of Lemma \ref{lm:cluster-expansion}]
    First by the uniform estimate, $|\varphi(\theta_p,Q'_p)| \leq 10^{4-6d}$, for any cluster $K$, with $|K|=m$, we have,
    \begin{equs}
        \bigg|\int  f(\theta,Q')\prod_{p \in K}\varphi(\theta_p,Q_p')\prod_{e \in E_K} d\nu_e(\theta_e) \bigg| \leq 10^{(4-6d)|K|}\|f\|_{L^{\infty}}.
    \end{equs}
    Similarly, by iterating Lemma \ref{lm:Comparison-of-partition-functions}, we see that 
    \begin{equs}            \frac{Z_{\Lambda \backslash \overline{K \cup \Lambda_f}}(Q')}{Z_{\Lambda}(Q')} \leq (1-10^{4-6d)})^{-|K|-|\Lambda_f|} \leq 2^{|K|+|\Lambda_f|}.
    \end{equs}
    Inputting the previous two estimates,
    \begin{equs}
        ~&\sum_{K \in \mc{C}(f):|K|>R}\bigg|\int  f(\theta,Q')\prod_{p \in K}\varphi(\theta_p,Q_p')\prod_{e \in E_K} d\nu_e(\theta_e) \frac{Z_{\Lambda \backslash \overline{K \cup \Lambda_f}}(Q')}{Z_{\Lambda}(Q')}\bigg| \\
        &=\sum_{m=R+1}^{\infty}\sum_{K \in \mc{C}(f):|K|=m}\bigg|\int  f(\theta,Q')\prod_{p \in K}\varphi(\theta_p,Q_p')\prod_{e \in E_K} d\nu_e(\theta_e) \frac{Z_{\Lambda \backslash \overline{K \cup \Lambda_f}}(Q')}{Z_{\Lambda}(Q')}\bigg|\\
        &\leq \sum_{m=R+1}^{\infty} 2^{|\Lambda_f|} (2 \cdot 10^{4-6d)})^m |\{K \in \mc{C}(f): |K|=m\}|\\
        &\leq 2^{|\Lambda_f|}\cdot e^{2d|\Lambda_f|}\sum_{m=R+1}^{\infty} (40^d \cdot 2 \cdot 10^{4-6d})^{m} \leq 2^{|\Lambda_f|}\cdot e^{2d|\Lambda_f|}\sum_{m=R+1}^{\infty} 10^{-6m(1-d)} \\
        & \leq  C(d,|\Lambda_f|) 10^{-3dR},
    \end{equs}
where Lemma \ref{lm:Combo-Lemma} was applied in the third to last line above, and that fact that $2 \cdot 40^d \leq 10^{2d}$, and $10^{6(1-d)m} \leq 10^{-3dm}$ for $d \geq 2$.
\end{proof}

For the next step in proving Lemma \ref{lm-Sensitivity-to-boundary-conditions} we first need some additional notation.

Fix $\Lambda' \subset \Lambda$. Then for a set of boundary conditions $(\theta_{e}')_{e \in (\Lambda')^c}$, and a cluster $K \subset \Lambda$, we set $Z_{\Lambda \backslash K,\theta'}$ to be the partition function,
    \begin{equs}
        Z_{\Lambda \backslash K,\theta'}=\int \prod_{p \in \mc{P}_{\Lambda \backslash K}^+}(1+\varphi(\theta_p,Q'_e))\prod_{e \in E_{(\Lambda')^c}^+} \delta(\theta_e-\theta_e') \prod_{e \in E_{\Lambda}^+} d\theta_e.
    \end{equs}

Before stating the next lemma, we remark that all the lemmas of the subsection up to this point apply equally well to the setting with boundary conditions described in the last paragraph.
    
\begin{lemma}\label{lm:difference-of-partition-function-ratio}
    For $\beta < \beta^*(d)$, and given the setup of the last paragraph, suppose there are two distinct boundary conditions $(\theta_e')_{e \in (\Lambda')^c}$ and $(\theta_e'')_{e \in (\Lambda')^c}$ and a cluster $K \subset \Lambda'$, then,
    \begin{equs}
        \bigg|\frac{Z_{\Lambda \backslash K, \theta'}}{Z_{\Lambda, \theta'}}-\frac{Z_{\Lambda \backslash K, \theta''}}{Z_{\Lambda, \theta''}}\bigg| \leq C_0(d)\cdot 2^{|K|}\cdot 10^{-2d\cdot d(K,(\Lambda')^c))}.
    \end{equs}
\end{lemma}

\begin{proof}
    The idea of the argument is to apply a double induction argument on the parameters $d(K,(\Lambda')^c)$ and $|\mc{P}_{\Lambda' \backslash K}^+|$. For the base case when $d(K,(\Lambda')^c)=0$ (note that we must have $d(K,(\Lambda')^c)=0$ when $|\mc{P}_{\Lambda' \backslash K}^+|=0$), we trivially have
    \begin{equs}
        ~&\bigg|\frac{Z_{\Lambda \backslash K, \theta'}}{Z_{\Lambda, \theta'}}-\frac{Z_{\Lambda \backslash K, \theta''}}{Z_{\Lambda, \theta''}}\bigg| \\
        &\leq \bigg|\frac{Z_{\Lambda \backslash K, \theta'}}{Z_{\Lambda, \theta'}}\bigg|+ \bigg|\frac{Z_{\Lambda \backslash K, \theta''}}{Z_{\Lambda, \theta''}}\bigg| \\
        &\leq C(d)2^{|K|}
    \end{equs}
    by applying Lemma \ref{lm:Comparison-of-partition-functions}.
    
    For the inductive step suppose $K=\{p\}$ is a single plaquette, and expand the partition functions in a cluster expansion as follows,
    \begin{equs}
        ~&\frac{Z_{\Lambda, \theta'}}{Z_{\Lambda \backslash \{p\}, \theta'}}-\frac{Z_{\Lambda, \theta''}}{Z_{\Lambda \backslash \{p\}, \theta''}}\\ 
        &=\sum_{J \subset \mc{C}(\{p\}): |J| <d(\{p\},(\Lambda')^c)} \int \prod_{p' \in J} \varphi(\theta_{p'} ,Q'_{p'}) \prod_{e \in E_J^+} d\theta_e \bigg(\frac{Z_{\Lambda \backslash (J \cup \{p\}), \theta'}}{Z_{\Lambda \backslash \{p\}, \theta'}}-\frac{Z_{\Lambda \backslash (J \cup \{p\}), \theta''}}{Z_{\Lambda \backslash \{p\}, \theta''}}\bigg)\\
        &+ \sum_{J \subset \mc{C}(\{p\}): |J| \geq d(\{p\},(\Lambda')^c)} \int \prod_{p' \in J} \varphi(\theta_{p'}, Q'_{p'}) \prod_{e \in E_J^+ \cap E_{(\Lambda')^c}^+}\delta(\theta_e-\theta_e') \prod_{e \in E_J^+} d\theta_e \frac{Z_{\Lambda \backslash (J \cup \{p\}), \theta'}}{Z_{\Lambda \backslash \{p\}, \theta'}}\\
        &-\sum_{J \subset \mc{C}(\{p\}): |J| \geq d(\{p\},(\Lambda')^c)} \int \prod_{p' \in J} \varphi(\theta_{p'}, Q'_{p'}) \prod_{e \in E_J^+ \cap E_{(\Lambda')^c}^+}\delta(\theta_e-\theta_e'') \prod_{e \in E_J^+} d\theta_e \frac{Z_{\Lambda \backslash (J \cup \{p\}), \theta''}}{Z_{\Lambda \backslash \{p\}, \theta''}}\\
        &=: S_1+S_2-S_3.
    \end{equs}
    To estimate $S_1$, recall $|\varphi(\theta_{p'} Q'_{p'})| \leq 10^{4-6d}$, so applying this $L^{\infty}$ bound together with the inductive hypothesis, and the combinatorial bound Lemma \ref{lm:Combo-Lemma}, we have,
    \begin{equs}
        |S_1| &\leq \sum_{J \subset \mc{C}(\{p\}): |J| <d(\{p\},(\Lambda')^c)} \int \prod_{p' \in J} |\varphi(\theta_{p'} ,Q'_{p'})| \prod_{e \in E_J^+} d\theta_e \bigg|\frac{Z_{\Lambda \backslash (J \cup \{p\}), \theta'}}{Z_{\Lambda \backslash \{p\}, \theta'}}-\frac{Z_{\Lambda \backslash (J \cup \{p\}), \theta''}}{Z_{\Lambda \backslash \{p\}, \theta''}}\bigg|\\
        &\leq \sum_{m=0}^{d(p,(\Lambda')^c)}C_0(d)40^{md} \cdot 10^{m(4-6d)} \cdot 2^m \cdot 10^{-2d(d(p,(\Lambda')^c)-m)}\\
        &\leq\sum_{m=0}^{d(p,(\Lambda')^c)}C_0(d) 2^{-m} 10^{-2d \cdot d(p,(\Lambda')^c)} \\
        &\leq C(d) \cdot 10^{-2d \cdot d(p,(\Lambda')^c)},
    \end{equs}
    where in the second to last line, we used the fact that $40^{md} \cdot 10^{m(4-6d)} \cdot 2^m \cdot 10^{2 md}\leq 2^{-m}$ for $d \geq 2$, and when applying the inductive bound above, we used the fact that $d(J,(\Lambda')^c) \leq d(p,(\Lambda')^c)-|J|$ as $J$ is a connected cluster of $p$.
     Similarly we can applying the bound $|\varphi(\theta_{p'} Q'_{p'})| \leq 10^{4-6d}\leq 10^{-4d}$ (for $d\geq 2$) together with the partition comparison Lemma \ref{lm:Comparison-of-partition-functions},
    \begin{equs}
        |S_2|+|S_3| \leq C(d)\sum_{m=d(p,(\Lambda')^c)}^{\infty} 40^{md}\cdot 10^{-4md} \cdot 2^{m} \leq C(d)10^{-2d(p,(\Lambda')^c)}.
    \end{equs}
     So finally setting $x=\frac{Z_{\Lambda \backslash \{p\}, \theta'}}{Z_{\Lambda, \theta'}}, y=\frac{Z_{\Lambda \backslash \{p\}, \theta''}}{Z_{\Lambda, \theta''}}$ we have shown that $|\frac{1}{x}-\frac{1}{y}|\leq C(d)10^{-2 d \cdot d(p,(\Lambda')^c)}$, and by Lemma \ref{lm:Comparison-of-partition-functions}, $|xy|\leq 4$, thus 
     \begin{equs}
         \bigg|\frac{Z_{\Lambda \backslash \{p\}, \theta'}}{Z_{\Lambda, \theta'}}-\frac{Z_{\Lambda \backslash \{p\}, \theta''}}{Z_{\Lambda, \theta''}}\bigg| =|x-y|=|xy||\frac{1}{x}-\frac{1}{y}|\leq C(d)10^{-2 d \cdot d(p,(\Lambda')^c)}.
     \end{equs}

     Finally, to get the general case, we write $K=\{p_1,\dots,p_n\}$ and let $x_k:=\frac{Z_{\Lambda' \backslash \{p_1,\dots,p_{k-1}\}, \theta'}}{Z_{\Lambda'\backslash \{p_1,\dots,p_k\}, \theta'}}$ and $y_k:=\frac{Z_{\Lambda' \backslash \{p_1,\dots,p_{k-1}\}, \theta''}}{Z_{\Lambda'\backslash \{p_1,\dots,p_k\}, \theta''}}$,
    \begin{equs}\label{eq:many-plaquette-difference}
        ~&\frac{Z_{\Lambda' \backslash K, \theta'}}{Z_{\Lambda', \theta'}}-\frac{Z_{\Lambda' \backslash K, \theta''}}{Z_{\Lambda', \theta''}} \\
        &= x_1x_2\cdots x_n - y_1y_2\cdots y_n =(x_1-y_1)x_2\cdots x_n + y_1(x_2-y_2)x_3 \cdots x_n \\
        &+ \dots +y_1 \cdots y_{n-1}(x_n-y_n),
    \end{equs}
    and by Lemma \ref{lm:Comparison-of-partition-functions}, $\max(|x_1|,\dots,|x_n|,|y_1|,\dots,|y_n|)\leq 2$, while by the previously discussed $K=\{p\}$ one plaquette case, $\max_i |x_i-y_i| \leq C(d)10^{-2d \cdot d(p_i,(\Lambda')^c)}$. The result now follows by plugging these bounds into \eqref{eq:many-plaquette-difference}. 
\end{proof}

We are now ready to establish the sensitivity of conditional expectations to boundary conditions.

\begin{proof}[Proof of Lemma \ref{lm-Sensitivity-to-boundary-conditions}]
    We first use the cluster expansion to expand,
    \begin{equs}
        ~&\E[f|Q=Q',(\theta_e)_{e \in E_{\Lambda^c}^+}=(\theta'_e)_{e \in E_{\Lambda^c}^+}]-\E[f|Q=Q',(\theta_e)_{e \in E_{\Lambda^c}^+}=(\theta''_e)_{e \in E_{\Lambda^c}^+}] \\
        &=\sum_{J \subset \mc{C}(\Lambda_f): |J| <d(\Lambda_f,(\Lambda')^c)} \int f(\theta,Q')\prod_{p \in J} \varphi(\theta_{p} ,Q'_{p}) \prod_{e \in E_J^+} d\theta_e \bigg(\frac{Z_{\Lambda \backslash (J \cup \Lambda_f), \theta'}}{Z_{\Lambda , \theta'}}-\frac{Z_{\Lambda \backslash (J \cup \Lambda_f), \theta''}}{Z_{\Lambda , \theta''}}\bigg)\\
        &+ \sum_{J \subset \mc{C}(\Lambda_f): |J| \geq d(\Lambda_f,(\Lambda')^c)} \int f(\theta,Q')\prod_{p \in J} \varphi(\theta_{p} ,Q'_{p}) \prod_{e \in E_J^+ \cap E_{(\Lambda')^c}^+}\delta(\theta_e-\theta_e') \prod_{e \in E_J^+} d\theta_e \frac{Z_{\Lambda \backslash (J \cup \Lambda_f), \theta'}}{Z_{\Lambda , \theta'}}\\
        &-\sum_{J \subset \mc{C}(\Lambda_f): |J| \geq d(\Lambda_f,(\Lambda')^c)} \int f(\theta,Q')\prod_{p \in J} \varphi(\theta_{p} ,Q'_{p}) \prod_{e \in E_J^+ \cap E_{(\Lambda')^c}^+}\delta(\theta_e-\theta_e'') \prod_{e \in E_J^+} d\theta_e \frac{Z_{\Lambda \backslash (J \cup \Lambda_f), \theta''}}{Z_{\Lambda, \theta''}}\\
        =: S_1+S_2-S_3.
    \end{equs}
    The required estimate for $|S_1|$ is similar to that in the proof of Lemma \ref{lm:difference-of-partition-function-ratio}, except we estimate 
    \begin{equs}
        \bigg|\frac{Z_{\Lambda \backslash (J \cup \Lambda_f), \theta'}}{Z_{\Lambda , \theta'}}-\frac{Z_{\Lambda \backslash (J \cup \Lambda_f), \theta''}}{Z_{\Lambda , \theta''}}\bigg|
    \end{equs}
    directly from Lemma \ref{lm:difference-of-partition-function-ratio} instead of via induction. The estimate of  $|S_2|+|S_3|$ follows exactly as in the analogous bound for the proof of Lemma \ref{lm:difference-of-partition-function-ratio}.
\end{proof}

\subsection{Step 2: Localizing conditional observables}\label{subsection:Step2}

The proof of Proposition \ref{prop:Localizing conditional observables}, starts by computing the derivative of the conditional expectation of local observable $f$ with respect to edge variables not belonging to $\Lambda_f$.

\begin{lemma}\label{lm:deriviative-of-conditional-expectation}
    Let $f$ be a local observable, $Q' \in \mrm{SU}(N)^{E_{\Lambda}^+}$ a fixed field, and $\tilde{X}$ be a tangent vector in $\mrm{T}_{Q'}(\mrm{SU}(N)^{E_{\Lambda}^+})$ of the form $(X_eQ_e)_{e \in E_{\Lambda}^+}$ for $X_e \in \mf{su}(N)$ such that $X_e = 0$ for $e \in \Lambda_f$, then,
    \begin{equs}
        \tilde{X}(\E[f|Q=Q']) = \sum_{e \in E_{\Lambda^+}} \sum_{\substack{p \in \mc{P}_{\Lambda}^+:\\p \succ e}} \mrm{sgn}(e,p) (\mrm{Cov}(f, \mrm{Re}(e^{i\frac{\theta_p}{N}}\mrm{Tr}(X_e Q_{p }))| Q=Q'),
    \end{equs}
    where $p \succ e$ denotes a plaquette containing $e$ or $e^{-1}$ as the first edge appearing in $p$.
\end{lemma}

\begin{proof}
    The proof is a direct calculation. First letting $Z'(Q'):= \int \exp(S_{\mc{U}}(\theta,Q))d\theta$, then
    \begin{equs}
        \tilde{X}(\tilde{Z}(Q')) &=\frac{d}{dt}|_{t=0} Z'((e^{tX_e}Q'_e)_{e \in E_{\Lambda}^+}]\\
        & = \sum_{e \in E_{\Lambda}^+} \int \sum_{p \succ e} \mrm{sgn}(e,p)\mrm{Re}(e^{i\frac{\theta_p}{N}}\mrm{Tr}(X_eQ'_p)) \exp(S_{\mc{U}}(\theta,Q'))d\theta.
    \end{equs}
    So by an additional similar computation,
    \begin{equs}
        \tilde{X}(\E[f|Q=Q'])&=\frac{d}{dt}|_{t=0} \E[f|Q_e=e^{tX_e}Q'_e, \text{ for each } e \in E_{\Lambda}^+]\\
        & =\sum_{e \in E_{\Lambda}^+} \int f(\theta,Q)\sum_{p \succ e} \mrm{sgn}(e,p)\mrm{Re}(e^{i\frac{\theta_p}{N}}\mrm{Tr}(X_eQ'_p)) \exp(S_{\mc{U}}(\theta,Q'))d\theta Z'(Q')^{-1}\\
        &-Z'(Q')^{-2}\sum_{e \in E_{\Lambda}^+} \int \sum_{p \succ e} \mrm{sgn}(e,p)\mrm{Re}(e^{i\frac{\theta_p}{N}}\mrm{Tr}(X_eQ'_p)) \exp(S_{\mc{U}}(\theta,Q'))d\theta \E[f|Q=Q']\\
        &=\sum_{e \in E_{\Lambda^+}} \sum_{p \succ e} \mrm{sgn}(e,p) (\mrm{Cov}(f, \mrm{Re}(e^{i\frac{\theta_p}{N}}\mrm{Tr}(X_e Q_{p }))| Q=Q').
    \end{equs}
\end{proof}

As a rather direct corollary of the formula from the last lemma and the conditional mass gap of Proposition \ref{prop:conditional-mass-gap}, we have an exponential decay bound on derivatives of the conditional expectation of $f$ with respect to edges far from $\Lambda_f$.

\begin{cor}\label{cor:derivative-sensitivity-to-Q}
    Let $f$ be a local observable, $Q' \in \mrm{SU}(N)^{E_{\Lambda}^+}$ a fixed field, and $\tilde{X}=(X_eQ_e)_{e \in E_{\Lambda}^+}$ be a tangent vector in $\mrm{T}_{Q'}(\mrm{SU}(N)^{E_{\Lambda}^+})$ such that $X_e = 0$ for all $e$ such that $d(e,\Lambda_f)<R$. Then there are constants $C_1=C_1(d,N,|\Lambda_f|),C_2=C_2(d,N)$ such that,
    \begin{equs}
        |\E[\tilde{X}(\E[f|Q=Q'])]| \leq C_1 e^{-C_2 R} \sup_{e} |X_e| \|f\|_{L^{\infty}}.
    \end{equs}
\end{cor}

\begin{proof}
    The strategy is to apply the previous lemma along with the conditional mass gap, Proposition \ref{prop:conditional-mass-gap}. For simplicity of notation, let $A_k=\{e \in E_{\Lambda^+}: R+k \leq d(e,\Lambda_f)< R+k+1 \}$ (the set of edges in the annulus at scale $R+k$). Then we have,
    \begin{equs}
        |\E[X(\E[f|Q=Q'])]| &\leq \sum_{e \in E_{\Lambda^+}: d(e,\Lambda_f) \geq R} \sum_{p \succ e} |(\E[\mrm{Cov}(f, \mrm{Re}(e^{i\frac{\theta_p}{N}}\mrm{Tr}(X_e Q_{p}))| Q=Q')]|\\
        &\leq \sum_{e \in E_{\Lambda^+}: d(e,\Lambda_f) \geq R} C_1\|\mrm{Re}(e^{i\frac{\theta_p}{N}}\mrm{Tr}(X_e Q_{p}))\|_{L^{\infty}} \|f\|_{L^{\infty}} e^{-C_2 d(e,\Lambda_f)} \\
        &\leq \sum_{k=0}^{\infty} C_1\sup_{e \in A_k}|X_e| \|f\|_{L^{\infty}} |A_k|e^{-C_2 (R+k))}\\
        &\leq \sum_{k=0}^{\infty} C_1'\sup_{e \in A_k}|X_e| \|f\|_{L^{\infty}} (R+k+1)^d e^{-C_2 (R+k))}\\
        &\leq C_1'' e^{-C_2' R}\sup_{e \in E_{\Lambda}^+} |X_e|\|f\|_{L^{\infty}},
    \end{equs}
    where the final constants $C_1''$ and $C_2'$ have the dependencies, $C_1''=C_1''(d,N,|\Lambda_f|)$, $C_2=C_2(d,N)$, and we used the fact that $\|\mrm{Re}(e^{i\frac{\theta_p}{N}}\mrm{Tr}(X_e Q_{p }))\|_{L^{\infty}} \leq \sqrt{N}|X_e|$ by Cauchy-Schwarz.
\end{proof}

\begin{proof}[Proof of Proposition \ref{prop:Localizing conditional observables}]
    Let $(Q'_e)_{e \in E_{\Lambda}^+}$ and $(Q''_e)_{e: d(e,\Lambda_f)\geq R}$ be two deterministic configurations, and consider a path $\gamma:[0,1] \to \mrm{SU}(N)^{\{e \in E_{\Lambda}^+: d(e,\Lambda_f)\geq R\}}$ such that $\gamma(0)=(Q'_e)_{e: d(e,\Lambda_f)\geq R}$ and $\gamma(1)=Q''$. Now since the diameter of $\mrm{SU}(N)$ is finite, (i.e. $\mrm{diam}(\mrm{SU}(N))\leq c_N$), we may choose $\gamma$ so that $|\gamma'(t)_e| \leq c_N$ for each $e \in E_{\Lambda}^+$. Thus by the fundamental theorem of calculus, and Corollary \ref{cor:derivative-sensitivity-to-Q},
    \begin{equs}
        |&\E[f| Q=Q']-\E[f|Q_e=Q_e' \text{ for } d(e, \Lambda_f)<R, \text{ and } Q_e=Q_e'' \text{ for } d(e, \Lambda_f)\geq R]| \\
        &=\bigg|\int_{0}^{1} \gamma'(t)(\E[f|Q_e=Q_e' \text{ for } d(e, \Lambda_f)<R, \text{ and } Q_e=\gamma(t)_e \text{ for } d(e, \Lambda_f)\geq R])dt \bigg|\\
        & \leq C_1 e^{-C_2 R}\sup_{t,e} |\gamma'(t)_e| \leq C_1' e^{-C_2 R},
    \end{equs}
    where we still have $C_1'=C_1'(d,N,|\Lambda_f|)$, $C_2=C_2(d)$. Now \eqref{eq:one-func-sensativity} follows by integrating with respect to the $Q''$ field. 
    
    Lastly Recall the definitions $\Lambda_{f,g^c}:=\{e:d(e,\Lambda_f) \leq \frac{1}{3}d(\Lambda_f,\Lambda_g)\}$ and $\Lambda_{f^c,g}:=\{e:d(e,\Lambda_g) \leq \frac{1}{3}d(\Lambda_f,\Lambda_g)\}$. By the same argument as the preceding paragraph, 
    \begin{equs}
        |\E[f|Q=Q']- \E[f|(Q_e)_{e \in \Lambda_{f,g^c}}]| \leq C_1 \|f\|_{L^{\infty}} e^{-C_2 d(\Lambda_f,\Lambda_g)},
    \end{equs}
    and
    \begin{equs}
        |\E[g|Q=Q']- \E[g|(Q_e)_{e \in \Lambda_{f^c,g}}]| \leq C_1'\|g\|_{L^{\infty}} e^{-C_2 d(\Lambda_f,\Lambda_g)},
    \end{equs}
    for constants $C_1=C_1(d,N,|\Lambda_f|)$ $C_1'=C_1'(d,N,|\Lambda_g|)$, and $C_2=C_2(d)$. Applying these estimates to the following decomposition establishes \eqref{eq:two-func-sensativity}. 
    \begin{equs}
        ~&|\E[f|Q=Q']\E[g|Q=Q'] -E[f|(Q_e)_{e \in \Lambda_{f,g^c}}]\E[g|(Q_e)_{e \in \Lambda_{f^c,g}}] |\\
        &\leq |\E[f|Q=Q']- \E[f|(Q_e)_{e \in \Lambda_{f,g^c}}]| |\E[g|Q=Q']|\\
        &+|\E[f|(Q_e)_{e \in \Lambda_{f,g^c}}]||\E[g|Q=Q']- \E[g|(Q_e)_{e \in \Lambda_{f^c,g}}]|\\
        &\leq (C_1+C_1')\|f\|_{L^{\infty}}\|g\|_{L^{\infty}}e^{-C_2 d(\Lambda_f,\Lambda_g)}.
    \end{equs}
\end{proof}

\subsection{Step 3: Mass gap for $\mrm{SU}(N)$ Marginal Field}\label{subsection:Step3}

The goal of this subsection is to prove Proposition \ref{prop:averaged-out-mass-gap}. Before diving into the proof let us recall that the marginal distribution of the $\mrm{SU}(N)$ can be expressed in the following form
\begin{equs}
    d\tilde{\mu}_{\mrm{SU}(N),\Lambda_L,\beta}(Q):=\frac{1}{\tilde{Z}_{\mrm{SU}(N),\Lambda_L,\beta}} \exp(\tilde{S}(Q)) dQ, 
\end{equs}
with
\begin{equs}
    \tilde{S}(Q):=\log\bigg( \int \exp(S_{\mc{U}}(\theta,Q))d\theta\bigg).
\end{equs}

The mass gap for this measure will be established through a similar Langevin dynamic based approach as in\cite{shen2023stochastic}. Due to the similarity with \cite{shen2023stochastic,shen2022new}, we focus on highlighting the differences, and on occasion we refer the reader back to these papers for more detail. The Langevin dynamics can informally be expressed as,
\begin{equs}
    dQ = \nabla \tilde{S}(Q) dt+ \sqrt{2} d\mathfrak{B},
\end{equs}
where $\mathfrak{B}$ denotes $\mrm{SU}(N)^{E_{\Lambda}^+}$-valued standard Brownian motion.

More explicitly we can write the dynamics above as the following system of stochastic differential equations.
\begin{equs}\label{eq:Langevin-Dynaimics-marginal}
    dQ_e = \nabla_e \tilde{S}(Q) dt + c_{\mathfrak{su}(N)} Q_e dt + \sqrt{2} dB_e Q_e,
\end{equs}
for each edge $e \in E_{\Lambda}^+$, and where each $B_e$ is an independent $\mf{su}(N)$-valued Brownian motion and $c_{\mathfrak{su}(N)}=\frac{1-N^2}{N}$ is the Casimir element for $\mathfrak{su}(N)$.

The existence and uniqueness for this system is quite standard.

\begin{lemma}[Global Well-Posedness]
    For any initial data $Q(0)=(Q_e(0))_{e \in E_{\Lambda}^+} \in \mrm{SU}(N)^{E_{\Lambda}^+}$, there exists a unique solution $Q=(Q_e)_{e \in E_{\Lambda}^+} \in C([0,\infty);\mrm{SU}(N)^{E_{\Lambda}^+})$ to \eqref{eq:Langevin-Dynaimics-marginal} a.s.
\end{lemma}
\begin{proof}
    Note that,
    \begin{equs}
        \exp(-N^2\beta |\mc{P}_{\Lambda}^+| )\leq \int \exp(S_{\mc{U}}(\theta,Q))d\theta \leq \exp(N^2\beta |\mc{P}_{\Lambda}^+| )
    \end{equs}
    and $\int \exp(S_{\mc{U}}(\theta,Q))d\theta$ is a smooth function of $Q \in \mrm{SU}(N)^{E_{\Lambda}^+}$, thus $\tilde{S}(Q)=\log \int \exp(S_{\mc{U}}(\theta,Q))d\theta$ is also smooth, and hence locally Lipschitz. The rest of the proof proceeds via standard arguments in SDE theory, as the coefficients of the SDE are locally Lipschitz and  $\mrm{SU}(N)^{E_{\Lambda}^+}$ is compact. The argument is sketched in \cite[Lemma 3.2]{shen2022new} and is based on arguments given in \cite{HsuStochasticGeometry,Levy2011a}.
\end{proof}

The main utility of the Langevin dynamics we study \eqref{eq:Langevin-Dynaimics-marginal}, is that we will later see that it mixes exponentially fast to the marginal $\mrm{SU}(N)$ field. Before proving this mixing property, we first establish that $\tilde{\mu}_{\mc{U}_L,\beta}$ is an invariant measure of the dynamics.

\begin{lemma}[Invariant Measure]
    The measure $\tilde{\mu}_{\mc{U}_L,\beta}$ defined in equation \eqref{eq:SU(N)-marginal-dist} is an invariant measure for the SDE system \eqref{eq:Langevin-Dynaimics-marginal}.
\end{lemma}

\begin{proof}
    The proof follows by integration by parts as in \cite[Lemma 3.3]{shen2023stochastic}.
\end{proof}

By global well-posedness of the SDE \eqref{eq:Langevin-Dynaimics-marginal}, the solutions form a Markov process in $\mrm{SU}(N)^{E_{\Lambda}^+}$.
We use $(P_t^{L})_{t\ge 0}=(P_t)_{t\ge 0}$ to denote the associated semigroup. For
$f\in C^\infty(\mrm{SU}(N)^{E_{\Lambda}^+})$,
\[
(P_t f)(Q_0)=\mathbb{E}\,f(Q(t,Q_0))
\qquad\text{for }Q_0\in \mrm{SU}(N)^{E_{\Lambda}^+},
\]
where $Q(t,Q_0)$ denotes the solution at time $t$ to \eqref{eq:Langevin-Dynaimics-marginal} starting from $Q_0\in \mrm{SU}(N)^{E_{\Lambda}^+}$.


Later on we will need the precise form of the gradient of the action appearing in the Langevin dynamics. This is the content of the next lemma.

\begin{lemma}\label{lm:gradient-Lipschitz-bound}
    \begin{equs}
        \nabla_e \tilde{S}(Q') =  N\beta \sum_{p \succ e} \mrm{sgn}(e,p)\E[\mbf{p}(e^{-i \frac{\theta_p}{N}} Q_p^*)Q_e|Q=Q'],
    \end{equs}
    where $\mbf{p}$ is the orthogonal projection of $\C^{N \times N}$ onto $\mf{su}(N)$, and recall that when we write $p \succ e$, we always require $e$ or $e^{-1}$ to be the first edge appearing in $p$.

    As a consequence for $\beta<\beta^*(d)=10^{-6d}$, and any edge $e$,
    \begin{equs}\label{eq:L^infty-gradient-bound}
        \|\nabla_e \tilde{S}(Q')\|_{L^{\infty}} \leq N^{3/2}.
    \end{equs}
\end{lemma}

\begin{proof}
    Let $v \in \mrm{T}_{Q}\mrm{SU}(N)^{E_{\Lambda}^+}$ which can be identified as $v=(X_eQ_e)_{e \in E_{\Lambda}^+}$ for $X_e \in \mf{su}(N)$. Then,
    \begin{equs}
        v(\tilde{S}(Q')) &= \E[N\beta \sum_{e \in E_{\Lambda}^+} \frac{d}{dt}|_{t=0} \sum_{p \succ e} \mrm{Re}( e^{i\frac{\theta_p}{N}} \mrm{Tr} (e^{\mrm{sgn}(e,p)X_e t} Q_p)|Q=Q'])\\
        &=\E[N\beta \sum_{e \in E_{\Lambda}^+}  \sum_{p \succ e} \mrm{sgn}(e,p)\mrm{Re} (e^{i\frac{\theta_p}{N}} \mrm{Tr} (X_e Q_p)|Q=Q'])\\
        &=\E[N\beta \sum_{e \in E_{\Lambda}^+}  \sum_{p \succ e} \mrm{sgn}(e,p)\langle X_e Q_e,e^{i \frac{\theta_p}{N}}Q_e^*Q_p\rangle|Q=Q']\\
        &=\sum_{e \in E_{\Lambda}^+}\langle v_e,  N\beta \sum_{p \succ e} \mrm{sgn}(e,p)\E[e^{-i \frac{\theta_p}{N}}Q_p^* Q_e|Q=Q']\rangle \\
        &= \langle v, N\beta \sum_{e \in E_{\Lambda}^+}  \sum_{p \succ e} \mrm{sgn}(e,p)\E[\mbf{p}(e^{-i \frac{\theta_p}{N}} Q_p^*)Q_e|Q=Q'] \rangle.
    \end{equs}
    From the definition of the gradient, $\nabla \tilde{S}(Q)$ is the unique vector in $\mrm{T}_{Q} \mrm{SU}(N)^{E_{\Lambda}^+}$ such that $\langle v,\nabla \tilde{S}(Q) \rangle = v(\tilde{S})(Q)$ for all tangent vectors $v \in \mrm{T}_{Q} \mrm{SU}(N)^{E_{\Lambda}^+}$. Thus the formula for $\nabla \tilde{S}(Q)$ follows.

    Next recall that the projection map $\mbf{p}$ simply maps a matrix $U \in \mrm{SU}(N)$ to $\frac{U-U^*}{2}-i(\mrm{Im}\mrm{Tr}U)I_N \in \mf{su}(N)$. So,
    \begin{equs}
        |\nabla_e \tilde{S}(Q)|
        &\leq N\beta\bigg(\sum_{p,p' \succ e}|\mrm{Re}\mrm{Tr}(\mbf{p} (e^{i\frac{\theta_p}{N}} Q_p) \mbf{p}(e^{-i\frac{\theta_{p'}}{N}} Q_{p'}^*))|\bigg)^{1/2}\\
        &\leq N\beta\bigg(\sum_{p,p' \succ e}\frac{1}{4}(|\mrm{Re}\mrm{Tr}(e^{i\frac{\theta_p}{N}} Q_p e^{-i\frac{\theta_{p'}}{N}} Q_{p'}^*)|+|\mrm{Re}\mrm{Tr}(e^{-i\frac{\theta_p}{N}} Q_p^* e^{-i\frac{\theta_{p'}}{N}} Q_{p'}^*)|\\
        &+|\mrm{Re}\mrm{Tr}(e^{i\frac{\theta_p}{N}} Q_p e^{i\frac{\theta_{p'}}{N}} Q_{p'})|+|\mrm{Re}\mrm{Tr}(e^{-i\frac{\theta_p}{N}} Q_p^* e^{i\frac{\theta_{p'}}{N}} Q_{p'})|\bigg)^{1/2},
    \end{equs}
    where we use the fact that $\mrm{Tr}(\mbf{p} (e^{i\frac{\theta_p}{N}} Q_p))=\mrm{Tr}(\mbf{p} (e^{i\frac{\theta_{p'}}{N}} Q_{p'}))=0$ (as these two matrices lie in $\mf{su}(N)$) to simplify the computation above.
    
    So finally for $\beta < \beta^*$
    each edge $e$, as there are at most $2(d-1)$ plaquettes containing $e$ as a side, we can use Cauchy-Schwarz to obtain the following bound
    \begin{equs}
        |\nabla_e \tilde{S}(Q)|
        \leq N\beta \cdot 2(d-1) \cdot N^{1/2} \leq 2(d-1)10^{-6d}N^{3/2} \leq N^{3/2},
    \end{equs}
    proving \eqref{eq:L^infty-gradient-bound}.
\end{proof}

Next we bound the Hessian of $\tilde{S}(Q)$.

\begin{lemma}\label{lm:Hessian-Bound}
    For any $\beta<\beta^*$ and any $v \in \mrm{T}_{Q}\mrm{SU}(N)^{E_{\Lambda}^+}$ for some $Q \in \mrm{SU}(N)^{E_{\Lambda}^+}$,
    \begin{equs}
        \mrm{Hess}(\tilde{S})(v,v) \leq C_d^* N \beta |v|^2,
    \end{equs}
    for a dimensional constant $C_d^*$.
\end{lemma}

\begin{proof}
    Throughout the proof, we identify $v \in \mrm{T}_Q\mrm{SU}(N)^{E_{\Lambda}^+}$ with $(X_eQ_e)_{e \in E_{\Lambda}^+}$ for $X_e \in \mathfrak{g}$. Next recall the formula for the Hessian \eqref{eq:Hessian-Formula},
    \begin{equs}
        \mrm{Hess}(\tilde{S})(v,v) = v(v(\tilde{S}))-(\nabla_v v)(\tilde{S})=v(v(\tilde{S})),
    \end{equs}
    where the second equality above follows from the explicit form of the Levi-Civita connection \eqref{eq:Levi-Civita-simplify}.
    
    Applying Lemma \ref{lm:gradient-Lipschitz-bound}, 
    \begin{equs}
        v(\tilde{S})=\langle v, \nabla \tilde{S}(Q)\rangle=\sum_{e \in E_{\Lambda}^+}  N\beta \sum_{p \succ e} \mrm{sgn}(e,p)\E[\mrm{Re}\mrm{Tr}(e^{i \frac{\theta_p}{N}}X_eQ_p) |Q=Q'].
    \end{equs}
    Next, a computation similar to that of Lemma \ref{lm:deriviative-of-conditional-expectation} shows that, 
    \begin{equs}
        v(v(\tilde{S}))&=N\beta \sum_{e,e' \in E_{\Lambda}^+}\sum_{\substack{p \in \mc{P}_{\Lambda}^+\\ p =ab}} \mrm{sgn}(e,p)\mrm{sgn}(e',p')\E[\mrm{Re}\mrm{Tr}(e^{i \frac{\theta_p}{N}} X_e Q_{a}X_{e'}Q_b)|Q=Q']\\
        &+\sum_{e,e' \in E_{\Lambda}^+}  (N\beta)^2 \sum_{\substack{p,p' \in \mc{P}_{\Lambda}^+ \\ p \succ e \\ p' \succ e'}} \mrm{sgn}(e,p)\mrm{sgn}(e',p') \mrm{Cov}(\mrm{Re}\mrm{Tr}(e^{i \frac{\theta_p}{N}}X_eQ_p),\mrm{Re}\mrm{Tr}(e^{i \frac{\theta_{p'}}{N}}X_{e'}Q_{p'}) |Q=Q')\\
        &=: I_1+I_2,
    \end{equs}
    where in the first line of the equation display above, $a$ and $b$ are possibly empty sequences of edges so that $p=ab$, and if $e \neq e'$, then $a$ is the sequence of edges in $p$ starting with $e$ or $e^{-1}$ and finishing immediately before $e'$ or $(e')^{-1}$, while if $e=e'$, then $b=p$ and $a$ is empty. For instance if $p=e_1e_2e_3e_4$ with $e_1=e$, and $e_3=e'$, then $a=e_1e_2$, $b=e_3e_4$, $Q_a=Q_{e_1}Q_{e_2}$, and $Q_b=Q_{e_3}Q_{e_4}$
    
    The sum defining $I_1$ can be bounded by $8N(d-1)\beta$ in the exact same manner as \cite[Lemma 4.1]{shen2023stochastic}. In particular note that $I_1=\E[F(\theta,Q)|Q=Q']$ for 
    \begin{equs}
        F(\theta,Q):=N\beta\sum_{e,e' \in E_{\Lambda}^+}\sum_{\substack{p \in \mc{P}_{\Lambda}^+\\ p \succ e,e'}} \mrm{sgn}(e,p)\mrm{sgn}(e',p')\mrm{Re}\mrm{Tr}(e^{i \frac{\theta_p}{N}} X_e Q_{a}X_{e'}Q_b),
    \end{equs}
    and by the Cauchy-Schwarz inequality with respect to the Frobenius norm, for every choice of $p,e,e'$
    \begin{equs}
        |\mrm{Re}\mrm{Tr}(e^{i \frac{\theta_p}{N}} X_e Q_{a}X_{e'}Q_b)| &\leq |X_eQ_{a}||X_{e'}Q_b|=|X_e||X_{e'}|\\
        &\leq \frac{1}{2}(|X_e|^2+|X_{e'}|^2).
    \end{equs}
    Thus summing over every valid choice of $p,e,e'$, we see that $|F(\theta,Q)| \leq 8(d-1)N\beta |v|^2$ where we picked up the combinatorial factor $8(d-1)$ due to the fact that every edge belongs to at most $2(d-1)$ positively oriented plaquettes. As a result $|I_1| \leq 8(d-1)N\beta |v|^2$

    Lastly to bound $I_2$, we first notice that by recentering the random variables in the covariance by the $\theta$ independent quantities $\mrm{Re}\mrm{Tr}(X_eQ_p)$ and $\mrm{Re}\mrm{Tr}(X_{e'}Q_{p'})$, we see that,
    \begin{equs}
        I_2=\sum_{e,e' \in E_{\Lambda}^+}  (N\beta)^2 \sum_{\substack{p \succ e \\ p' \succ e'}} \mrm{sgn}(e,p)\mrm{sgn}(e',p') \mrm{Cov}(\mrm{Re}\mrm{Tr}((e^{i \frac{\theta_p}{N}}-1)X_eQ_p),\mrm{Re}\mrm{Tr}((e^{i \frac{\theta_{p'}}{N}}-1)X_{e'}Q_{p'}) |Q=Q').
    \end{equs}
   Moreover since $|e^{i \frac{\theta_{p}}{N}}-1|,|e^{i \frac{\theta_{p'}}{N}}-1| \leq \frac{C}{N}$ for a constant $C$,
   \begin{equs}
       |\mrm{Re}\mrm{Tr}((e^{i \frac{\theta_p}{N}}-1)X_eQ_p)| &\leq \frac{C}{N}|\mrm{Tr}(X_eQ_p)|\\
       &\leq \frac{C}{\sqrt{N}} |X_eQ_p|=\frac{C}{\sqrt{N}} |X_e|,
   \end{equs}
   where we applied Cauchy-Schwarz to the inner product $\langle I_N,X_eQ_p \rangle$. Similarly we have $|\mrm{Re}\mrm{Tr}((e^{i \frac{\theta_{p'}}{N}}-1)X_{e'}Q_{p'})| \leq \frac{C}{\sqrt{N}}$. We are now ready to apply the conditional mass gap, Proposition \ref{prop:conditional-mass-gap} to deduce that,
   \begin{equs}
       ~&|\mrm{Cov}(\mrm{Re}\mrm{Tr}((e^{i \frac{\theta_p}{N}}-1)X_eQ_p),\mrm{Re}\mrm{Tr}((e^{i \frac{\theta_{p'}}{N}}-1)X_{e'}Q_{p'}) |Q=Q')| \leq \frac{C_1}{N} 10^{-2d \cdot d(p,p')} |X_e| |X_{e'}|\\
       &\leq \frac{C_1}{2N} 10^{-2d \cdot d(p,p')}(|X_e|^2+|X_{e'}|^2),
   \end{equs}
   for dimensional constants $C_1$ and $C_2$. So finally,
   
    \begin{equs}
        |I_2| \leq (N\beta)^2 \sum_{\substack{p \succ e \\ p' \succ e'}} \frac{C_1}{2N} 10^{-2d \cdot d(p,p')}(|X_e|^2+|X_{e'}|^2) \leq C(d)N \beta^2|v|^2 \leq C(d)N\beta |v|^2,
    \end{equs}
    for a new constant $C(d)$, and where we used $\beta <1$. 

    In summary we have shown that
    \begin{equs}
         ~&|\mrm{Hess}(\tilde{S})(v,v)| = |v(v(\tilde{S}))|\\
         &\leq |I_1| + |I_2| \leq C_d^* N\beta |v|^2,
    \end{equs}
    for some dimensional constant $C_d^*$.
\end{proof}

From now on we set $K_{\tilde{S}}:=\frac{N+2}{2}-1-C_d^* N \beta$ for the $C_d^*$ appearing in the statement of Lemma \ref{lm:Hessian-Bound}. Moreover, we can take $\tilde{\beta} = \min((3C_d^*)^{-1},\beta^*)$ so that $K_{\tilde{S}}$ is bounded below by a positive constant only depending on $N$ and $d$ whenever $\beta<\tilde{\beta}$. We are now in position to verify the Bakry-Émery criterion and thus establish exponential mixing of the Langevin dynamics \eqref{eq:Langevin-Dynaimics-marginal} as in \cite[Theorem 4.2]{shen2023stochastic}.

\begin{prop}\label{prop:Poincare-mixing}
Suppose $\beta < \tilde{\beta}$ ($\tilde{\beta}$ taken as in the preceding paragraph). The Langevin dynamics \eqref{eq:Langevin-Dynaimics-marginal} is exponentially ergodic in the sense that
\[
\bigl\| P_t^L f - \tilde{\mu}_{\mrm{SU}(N),\Lambda_L,\beta}(f) \bigr\|_{L^2(\tilde{\mu}_{\mrm{SU}(N),\Lambda_L,\beta})}
   \;\leq\; e^{-tK_{\tilde{S}}} \, \| f \|_{L^2(\tilde{\mu}_{\mrm{SU}(N),\Lambda_L,\beta})},
\]
where $ \tilde{\mu}_{\mrm{SU}(N),\Lambda_L,\beta}(f)$ is the expectation of $f$ with respect to $\tilde{\mu}_{\mrm{SU}(N),\Lambda_L,\beta}$. Moreover the invariant measure of $(P_t^L)_{t\ge0}$ is unique.
\end{prop}

\begin{remark}
    There are additional immediate consequences of the Bakry-Émery condition such as a log-Sobolev inequality, however we limit our attention only to the consequence necessary for the proof of a mass gap.
\end{remark}

\begin{proof}

As in the proof of \cite[Theorem 4.2, Remark 4.6] {shen2023stochastic}, the conditions of the theorem are equivalent to the Bakry-Émery condition: for every $v=XQ\in \mrm{T}_Q \mrm{SU}(N)^{E_{\Lambda}^+}$,
\begin{equation}\label{eq:Bakry-Emery}
\mathrm{Ric}(v,v) - \langle \nabla_v \nabla \tilde{S}, v\rangle \;\ge\; K_{\tilde{S}} |v|^2. 
\end{equation}
Here we recall that $|v|^2 = |X|^2$ and $\langle \nabla_v \nabla \tilde{S}, v\rangle
= \mathrm{Hess}(\tilde{S})(v,v)$. By \cite[(F.6)]{RMTBook}, for any tangent vector $u$ of $\mrm{SU}(N)$,
\[
\mathrm{Ric}(u,u) = \Bigl(\frac{(N+2)}{2}-1\Bigr)|u|^2.
\]

Since $\mathrm{Ric}(v,v)=\sum_e \mathrm{Ric}(v_e,v_e)$ and $|X|^2=\sum_e |X_e|^2$,
we have
\begin{equation}
\mathrm{Ric}(v,v) = \Bigl(\frac{(N+2)}{2}-1\Bigr)|X|^2. 
\end{equation}
As a result, using Lemma \ref{lm:Hessian-Bound} and the definition of $K_{\tilde{S}}$ we can immediately verify \eqref{eq:Bakry-Emery}.
\end{proof}

In the remainder of the section, it will be convenient for the calculations to consider an explicit choice of an
orthonormal basis of $\mathfrak{su}(N)$. This choice also appears in \cite[Proposition~E.15]{RMTBook} and \cite[Section~4.3]{shen2023stochastic}. Let $e_{kn}\in M_N$ for $k,n=1,\dots,N$ be the elementary matrices, namely the matrix whose
$(k,n)$-th entry is $1$ and all other entries are $0$.
For $1\le k<N$, and let
\[
D_k=\frac{i}{\sqrt{k+k^2}}\!\left(-k\,e_{k+1,k+1}+\sum_{j=1}^{k} e_{jj}\right).
\]
For $1\le k,n\le N$, let
\begin{equation}
E_{kn}=\frac{e_{kn}-e_{nk}}{\sqrt{2}},
\qquad
F_{kn}=\frac{i e_{kn}+i e_{nk}}{\sqrt{2}}.
\end{equation}
Then:
$\{D_k:1\le k<N\}\,\cup\,\{E_{kn},F_{kn}:1\le k<n\le N\}$ is an orthonormal basis of $\mathfrak{su}(N)$.

This then determines an orthonormal basis $\{v_e^{\,i}\}$ of $\mathfrak{su}(N)^{E_{\Lambda}^+}$,
which consists of right-invariant vector fields on $\mrm{SU}(N)^{E_{\Lambda}^+}$.

Now that we have verified exponential mixing of the Langevin dynamics, the next step towards the marginal mass gap is to prove that the dynamics is local in a quantitative sense which will be captured by Lemma \ref{lm:Generator-Commutator-Bound}. In preparation we first need to recall the estimate of Lemma \ref{lm:basis-commutator-bound}.

\begin{lemma}\label{lm:basis-commutator-bound}\cite{shen2023stochastic}[Lemma 4.9] 
    It holds for every $v_e^i$
    \begin{equs}
        \sum_j &|[v_e^i,v_e^j]f|^2 \leq \frac{9}{2} |\nabla_e f|^2 \hspace{3mm} \text{ for } \hspace{3mm} G= \mathrm{SU}(N),
    \end{equs}
\end{lemma}

\begin{lemma}\label{lm:Generator-Commutator-Bound}
Suppose $\beta < \tilde{\beta}(d)$. Let $\{v_e^{\,i}\}$ be the orthonormal basis given above. For every $f\in C^\infty(\mrm{SU}(N)^{E_{\Lambda}^+})$
and every $e\in E_{\Lambda}^+$, one has
\[
\bigl|[v_e^{\,i},\mathcal{L}]f(Q)\bigr|
\;\le\; \sum_{\bar e \in E_{\Lambda}^+} a_{e,\bar e}\,
\bigl|\nabla_{\bar e} f(Q)\bigr|,\qquad \forall\,Q\in \mrm{SU}(N)^{E_{\Lambda}^+},
\]
with $|a_{e,\bar e}| \leq C_1 10^{-2d \cdot d(e,\bar e)}$ for a universal constant $C=C(d,N)$.
\end{lemma}

\begin{proof}
The first few equation displays in the proof \cite[Lemma 4.10]{shen2023stochastic} apply verbatim to our setting, leading to the following identity,
\begin{align*}
[v_e^{\,i},\mathcal{L}]f
&= \sum_{\bar e \in E_{\Lambda}^+} \Bigl\langle \sum_j (v_e^{\,i}v_{\bar e}^{\,j}\tilde{S})\, v_{\bar e}^{\,j},\, \nabla_{\bar e}f \Bigr\rangle
+ \frac{3}{2}\sum_j (v_e^{\,j}\tilde{S}) \bigl\langle [v_e^{\,i},v_e^{\,j}],\, \nabla_e f \bigr\rangle \\
&\qquad + \frac12\sum_j v_e^{\,j}f \bigl\langle \nabla_e\tilde{S},\, [v_e^{\,i},v_e^{\,j}] \bigr\rangle
\; =:\; \sum_{k=1}^{3} I_k .
\end{align*}

We start by bounding $I_1$. By a similar calculation as in the proof of the Hessian bound Lemma \ref{lm:Hessian-Bound}, we have
$\lvert v_e^{\,i}v_{\bar e}^{\,j}\tilde{S}\rvert \le 8(d-1) N\beta \delta_{e \sim \bar e}+C_1 10^{-2d \cdot d(e,\bar e)}$ for $\bar e\neq e$. So substituting the basis expansion $\nabla_{\bar e} f = \sum_{k}(v_{\bar e}^k f)v_{\bar e}^k$, applying the Cauchy-Schwarz
inequality, and plugging in the second derivative bound of the previous sentence, we have,
\begin{align*}
\lvert I_1\rvert
&= \Bigl|\sum_{\bar e \in E_{\Lambda}^+}\sum_j (v_e^{\,i}v_{\bar e}^{\,j}\tilde{S})\, v_{\bar e}^{\,j}f \Bigr| \\
&\le \sum_{\bar e \in E_{\Lambda}^+} \Bigl(\sum_j \lvert v_e^{\,i}v_{\bar e}^{\,j}\tilde{S}\rvert^2 \Bigr)^{1/2}
\Bigl(\sum_j \lvert v_{\bar e}^{\,j}f\rvert^2 \Bigr)^{1/2} \\
&\le \sqrt{d(\mathfrak g)} \sum_{\bar e \in E_{\Lambda}^+} C_1 e^{-C_2 d(\bar e,e)}\lvert \nabla_{\bar e} f\rvert.
\end{align*}

Next we bound $I_2+I_3$. We first apply Cauchy-Schwarz inequality, then the definition of the gradient, and next Lemma \ref{lm:basis-commutator-bound} as follows,
\begin{align*}
\lvert I_2+I_3\rvert
&\le \frac{3}{2}\Bigl(\sum_j \lvert v_e^{\,j}\tilde{S}\rvert^2 \Bigr)^{1/2}
\Bigl(\sum_j \lvert\langle [v_e^{\,i},v_e^{\,j}],\, \nabla_e f\rangle\rvert^2 \Bigr)^{1/2}
+ \frac{1}{2}\Bigl(\sum_j \lvert v_e^{\,j}f\rvert^2 \Bigr)^{1/2}
\Bigl(\sum_j \lvert \langle [v_e^{\,i},v_e^{\,j}],\, \nabla_e\tilde{S}\rangle \rvert^2 \Bigr)^{1/2} \\
&=\frac{3}{2}\Bigl(\sum_j \lvert v_e^{\,j}\tilde{S}\rvert^2 \Bigr)^{1/2}
\Bigl(\sum_j \lvert [v_e^{\,i},v_e^{\,j}] f\rvert^2 \Bigr)^{1/2}
+ \frac{1}{2}\Bigl(\sum_j \lvert v_e^{\,j}f\rvert^2 \Bigr)^{1/2}
\Bigl(\sum_j \lvert [v_e^{\,i},v_e^{\,j}]\tilde{S}\rvert^2 \Bigr)^{1/2}\\
&\le 3\sqrt{2}\, \Bigl(\sum_j \lvert v_e^{\,j}\tilde{S}\rvert^2 \Bigr)^{1/2}
\Bigl(\sum_j \lvert v_e^{\,j}f\rvert^2 \Bigr)^{1/2}
= 3\sqrt{2}|\nabla_e\tilde{S}|\, |\nabla_e f| \\
&\le 3\sqrt{2}N^{3/2} |\nabla_e f|.
\end{align*}
 In the last line we used the bound on $\|\nabla_e \tilde{S}\|_{L^{\infty}}$ from equation \eqref{eq:L^infty-gradient-bound}.
\end{proof}

\noindent
We are now ready to prove the marginal distribution mass gap property of Proposition \ref{prop:averaged-out-mass-gap}.

\begin{proof}[Proof of Proposition \ref{prop:averaged-out-mass-gap}]

The argument is similar to that of \cite[Corollary 4.11]{shen2023stochastic}, with the main difference coming from the fact that the action $\tilde{S}$ considered here is not exactly local while the original Yang-Mills action is local. In the argument that follows, all expectations and covariances are taken with respect to $\tilde{\mu}_{\mrm{SU}(N),\Lambda_L,\beta}$ which we will denote by $\E_{\tilde{\mu}}$ and $\mrm{Cov}_{\tilde{\mu}}$ respectively. For simplicity of notation we will also drop the $L$ when we refer to the semigroup $P_t^L$ and the generator $\mc{L}_L$ For any two local observables $f,g \in C^{\infty}(\mrm{SU}(N)^{E_{\Lambda}^+})$ we have the following covariance decomposition coming from the law of total expectation,
\begin{equs}\label{eq:Covar-decomp}
    \mathrm{Cov}_{\tilde{\mu}}(f,g)
    &=\E_{\tilde{\mu}}[fg]-\E_{\tilde{\mu}}[f]E_{\tilde{\mu}}[g]
    =\E_{\tilde{\mu}}[P_t(fg)]-\E_{\tilde{\mu}}[P_t f]\E_{\tilde{\mu}}[P_t g] \\
    &=\E_{\tilde{\mu}}[P_t(fg)-P_t f\,P_t g]
    + \mathrm{Cov}_{\tilde{\mu}}(P_t f,P_t g),
\end{equs}
and by the Cauchy-Schwartz inequality,
\begin{equs}        
    |\mathrm{Cov}_{\tilde{\mu}}(P_t f,P_t g)|\leq\mathrm{Var}(P_t f)^{1/2}\,\mathrm{Var}(P_t g)^{1/2}
\end{equs}
Now recall that the Poincaré inequality Proposition \ref{prop:Poincare-mixing} is equivalent to
\[
    \mathrm{Var}(P_t f)\;\le\; e^{-2tK_{\tilde{S}}}\,\|f\|_{L^2(\tilde{\mu}_{\mc{U}_L,\beta})}^2,
\] 
for any $f \in C^{\infty}(\mrm{SU}(N)^{E_{\Lambda}^+})$Therefore by the Poincaré inequality, 
    is bounded by
\begin{equs}\label{eq:Mixing-Cov}
    |\mathrm{Cov}_{\tilde{\mu}}(P_t f,P_t g)|
    \;\le\; e^{-2tK_{\tilde{S}}}\,\|f\|_{L^2(\tilde{\mu}_{\mc{U}_L,\beta})}\|g\|_{L^2(\tilde{\mu}_{\mc{U}_L,\beta})}. 
\end{equs}

    Now we conside the term $P_t(fg)-P_t f P_t g$ in \eqref{eq:Covar-decomp} and we omit $L$
    for notation simplicity. The analysis begins by interpolating between time $0$ and $t$ as follows. In what follows we use the fact that $P_t$ and $\mathcal L$ commute
(see e.g.\ \cite[Chap.~I, Exercise 1.9]{DirichletFormsBook}), and that since $\mathcal{L}$ is a uniform elliptic operator with smooth coefficients,
    by Hörmander’s Theorem (cf.\ \cite[Theorem~2.3.3]{MalliavanCalc09}) $P_t f \in C^\infty(\mrm{SU}(N)^{E_{\Lambda}^+})$. By the fundamental theorem of calculus and Leibnitz rule,
\begin{equs}\label{eq:covar-first-part}
P_t(fg)-P_t f\,P_t g
&= \int_0^t \tfrac{d}{ds}\bigl[P_s(P_{t-s} f\, P_{t-s} g)\bigr]\,ds \\
&= \int_0^t \bigl[P_s \mathcal L(P_{t-s} f P_{t-s} g) - P_s(\mathcal L P_{t-s} f P_{t-s} g + P_{t-s} f \mathcal L P_{t-s} g)\bigr]\,ds \\
&= 2\sum_e \int_0^t P_s\langle \nabla_e P_{t-s} f, \nabla_e P_{t-s} g\rangle\,ds
= 2\sum_{e,i}\int_0^t P_s\bigl[(v_e^i P_{t-s} f)\cdot (v_e^i P_{t-s} g)\bigr]\,ds,
\end{equs}
where in the last line we used the fact that the first order derivative terms arising from $\mc{L}$ all cancel due to the Leibnitz rule, and moreover that for any $F,G \in C^{\infty}$, $\Delta_e (FG)-(\Delta_eF)G-F(\Delta_e G) = \langle\nabla_e F,\nabla_e G\rangle $. Finally in the last equality of the equation display above, we used the basis expansion $\nabla_e F=\sum_{j}(v_e^j F)v_e^j$, $\nabla_e G=\sum_{j}(v_e^j G)v_e^j$, along with the fact that $\{v_e^j\}_j$ is an orthonormal basis. 

Suppose for the moment that there is a constant $B$ (possibly depending on $N$ and $d$) such that for all $t<\frac{d(\Lambda_f,\Lambda_g)}{B}$ and $f \in C^\infty(\mrm{SU}(N)^{E_{\Lambda}^+})$, there are universal constants $C_1=C_1(|\Lambda_f|,N,d), C_2=C_2(N,d)>0$ one has
\begin{equs}\label{eq:basis-propogator-commutator}
\sup_{i=1,\dots,d(\mf{g})} \|v_e^i P_t f \|_{L^{\infty}}
\;\le\; C_1\,  e^{-2C_2 d(\Lambda_f,\Lambda_f)}\vertiii{f}_\infty,
\end{equs}
for all $e$ such that $d(e,\Lambda_f)\geq \frac{1}{3} d(\Lambda_f,\Lambda_g)$. If this were true, choosing $t=\frac{d(\Lambda_f,\Lambda_g)}{2B}$, by \eqref{eq:Mixing-Cov},
\begin{equs}\label{eq:final-cov-1}
    |\mrm{Cov}_{\tilde{\mu}}(P_t f,P_t g)| \leq e^{-\frac{K_{\tilde{S}}}{B}d(\Lambda_f),\Lambda_g)} \|f\|_{L^2(\tilde{\mu})}\|g\|_{L^2(\tilde{\mu})},
\end{equs}
while plugging \eqref{eq:basis-propogator-commutator} into \eqref{eq:covar-first-part}, and applying the conditional Jensen inequality,
\begin{equs}\label{eq:final-cov-2}
    \|P_t(fg)-P_t(g)P_t(g)\| &\leq \frac{d(\Lambda_f,\Lambda_g)}{B} \cdot d(\mf{g})\cdot \vertiii{f}_{\infty}\vertiii{g}_{\infty}\cdot C_1 e^{-C_2 d(\Lambda_f,\Lambda_g)}\\
    &\leq C_1' e^{-\frac{C_2}{2} d(\Lambda_f,\Lambda_g)}\vertiii{f}_{\infty}\vertiii{g}_{\infty},
\end{equs}
where we used the fact that by the triangle inequality, either $d(e,\Lambda_f)\leq \frac{1}{3} d(\Lambda_f,\Lambda_g)$ or $d(e,\Lambda_g)\leq \frac{1}{3} d(\Lambda_f,\Lambda_g)$, and the elementary inequality $x \leq \frac{2}{C_2}e^{\frac{C_2}{2}x}$ for $x,C_2 \geq 0$. Lastly plugging in the bounds \eqref{eq:final-cov-1} and \eqref{eq:final-cov-2} into the decomposition \eqref{eq:Covar-decomp} proves Proposition \eqref{prop:averaged-out-mass-gap} modulo \eqref{eq:basis-propogator-commutator}.



To verify the claimed bound \eqref{eq:basis-propogator-commutator}, we use the argument of
\cite[Theorem~8.2]{GuionnetLogSobolev} adapted to our setting. Using the fundamental theorem of calculus and the Leibnitz rule,
\begin{equs}
v_e^i P_t f - P_t v_e^i f
&= \int_0^t \frac{d}{ds} \bigl(P_{t-s} v_e^i P_s f\bigr)\,ds
= \int_0^t P_{t-s}[v_e^i,\mathcal L] P_s f\,ds. 
\end{equs}
Now recall Lemma \ref{lm:Generator-Commutator-Bound},
\[
\|[v_e^i,\mathcal L] P_s f\|_{L^{\infty}} \le \sum_{\bar e \in E_{\Lambda}^+} a_{e,\bar e}\|\nabla_{\bar e} P_s f\|_{L^{\infty}} ,
\]
where we recall that the coefficients $a_{e, \bar e} \leq C 10^{-2d \cdot d(e,e')}$ for some constants $C=C(N,d)$. Combining the previous two equation displays, we have
\begin{equs}\label{eq:integral-eq-mass-gap}
    \|v_e^i P_t f\|_{L^{\infty}}  \leq \| P_t v_e^i f\|_{L^{\infty}}+\sum_{j} \int_{0}^{t} \sum_{\bar e \in E_{\Lambda}^+} a_{e,\bar e}\|v_{\bar e}^j P_s f\|_{L^{\infty}} ds.
\end{equs}
Now for positive constants $T,\lambda>0$, let us define a norm,
\begin{equs}
    \|f\|_{T,\lambda}:=\sum_{e \in E_{\Lambda}^+} e^{\lambda d(e,\Lambda_f)} \sup_{s \in [0,T]}\sup_{i} |v_e^i P_s f\|_{L^{\infty}}.
\end{equs}
Set $\lambda =d\log(10)$. Then taking a weighted average of \eqref{eq:integral-eq-mass-gap} over all edges, we have,
\begin{equs}
    \|f\|_{t,\lambda} &\leq \sum_{e \in E_{\Lambda_f}^+}\sup_{s \in [0,t]}\sup_{i}\|P_s v_e^i f\|_{L^{\infty}}+\int_{0}^{t} \sum_{e,\bar e \in E_{\Lambda}^+} e^{\lambda d(e,\Lambda_f)}a_{e,\bar e}\sum_j\|v_{\bar e}^j P_s f\|_{L^{\infty}} ds\\
    &\leq C\vertiii{f}_{\infty}+\int_{0}^{t} \sum_{e,\bar e \in E_{\Lambda}^+} e^{\lambda d(e,\Lambda_f)}a_{e,\bar e}\sum_j\|v_{\bar e}^j P_s f\|_{L^{\infty}} ds\\
    &\leq C\vertiii{f}_{\infty}+C'\int_{0}^{t} d(\mf{g})\sum_{\bar e \in E_{\Lambda}^+} \sum_{e \in E_{\Lambda}^+} e^{\lambda d(e,\Lambda_f)-2d\log(10)d(e,\bar e)}\sup_j\|v_{\bar e}^j P_s f\|_{L^{\infty}} ds\\
    &\leq C\vertiii{f}_{\infty}+C''\int_{0}^{t} \sum_{\bar e \in E_{\Lambda}^+} \sum_{e \in E_{\Lambda}^+} e^{(\lambda-2d\log(10))d(e,\bar e)}\cdot e^{\lambda d(\bar e,\Lambda_f)}\sup_j\|v_{\bar e}^j P_s f\|_{L^{\infty}} ds\\
    &\leq C\vertiii{f}_{\infty}+C'''\int_{0}^{t} \sum_{\bar e \in E_{\Lambda}^+} e^{\lambda d(\bar e,\Lambda_f)}\sup_j\|v_{\bar e}^j P_s f\|_{L^{\infty}} ds\\
    &\leq C\vertiii{f}_{\infty}+C'''\int_{0}^{t} \|f\|_{s,\lambda} ds,
\end{equs}
where in the inequalities above $C,C',C'',C'''$ are all constant depending on $d$ and $N$, and we used the fact that $P_s v_e^i f \neq 0$ only if $e \in E_{\Lambda_f}^+$. We also applied the triangle inequality $d(e,\Lambda_f) \leq d(e,\bar e)+d(\bar e, \Lambda_f)$.

Applying Grönwall's inequality we have, with the same constants $C$ and $C'''$ as in the previous equation display,
\begin{equs}
    \|f\|_{t,\lambda} &\leq C\vertiii{f}_{\infty} +\int_{0}^{t} C\vertiii{f}_{\infty}C'''e^{C'''(t-s)}ds\\
    &\leq C(1+e^{C''' t})\vertiii{f}_{\infty}.
\end{equs}
So finally we make the choice $B=10C'''$, for $C'''$ as above, and  recall that $t<\frac{d(\Lambda_f,\Lambda_g)}{B}$. Since $\sup_i \|P_t v_e^i\|_{L^{\infty}} \leq e^{-\frac{\lambda}{3}d(\Lambda_f,\Lambda_g)}\|f\|_{t,\lambda}$ for all edges $e$ satisfying $d(e,\Lambda_f)\geq \frac{1}{3}d(\Lambda_f,\Lambda_g)$, we have that  \eqref{eq:basis-propogator-commutator} holds completing the proof.


\end{proof}

\section{Infinite Volume Limit}\label{section:infinite-volume}

In this section we emphasize the dependence on $L$, and sending $L\to \infty$ we extend $\Lambda_L$ to $\Lambda_{\infty} =\Z^d$. As in \cite[(3.6)]{shen2023stochastic}, we define the (squared) norm 
\begin{equs}
    \|Q\|^2:= \sum_{e \in E_{\Lambda_{\infty}}^+} \frac{1}{2^{|e|}} |Q_e|^2,
\end{equs}
on $M_N^{E_{\Lambda_{\infty}}^+}$, where $M_N$ is the space of $N \times N$ complex matrices, and where $|e|$ denotes the distance from $e$ to $0$ in $\Z^d$ and $|Q_e|$ denotes the usual Frobenius norm of the matrix $Q_e$. This norm induces a topology on $\mrm{U}(N)^{E_{\Lambda_\infty}^+} \sse M_N^{E_{\Lambda_\infty}^+}$. In the following, weak convergence of probability measures on $\mrm{U}(N)^{E_{\Lambda_\infty}^+}$ will be with respect to the Borel $\sigma$-algebra induced by this norm.

The lattice Yang-Mills theory on $\Lambda_L$ with inverse temperature $\beta$, $\mu_{\mrm{U}(N),\Lambda_L,\beta}$, can be extended to a measure on $\Lambda_{\infty}$ in an arbitrary way, for instance by assigning an independent Haar-distributed matrix to each edge outside $\Lambda_L$. In this section, we finish the proof of Theorem \ref{thm:U-N-mass-gap} by showing that for $\beta < \tilde{\beta}$ (the regime of $\beta$ up to which mass gap was proven), $\mu_{\mrm{U}(N),\Lambda_L,\beta}$ converges weakly to a measure $\mu_{\mrm{U}(N),\Lambda_{\infty},\beta}$ on $\mrm{U}(N)^{E_{\Lambda_{\infty}}^+}$ which we call the infinite volume limit. 

The first ingredient towards the existence of the infinite volume limit is a generalization of the mass gap statement for a mixed temperature model. The key observation here is that actually all of the proofs thus far both in this paper and in \cite{shen2023stochastic} apply equally well to the mixed temperature version of the model introduced in the mass gap lemma below.


\begin{lemma}[Mixed temperature mass gap]\label{lm:Mixed-Temp-Mass-Gap}
Let $\mbf{B} = (\beta_p)_{p \in \mc{P}_{\Lambda_L}^+}$ be a collection of non-negative numbers. Define the measure $\mu_{\mrm{U}(N),\Lambda_L,\mbf{B}}$ on $\mrm{U}(N)^{E_{\Lambda_L}^+}$ by
\begin{equs}
    d\mu_{\mrm{U}(N),\Lambda_L,\mbf{B}} = \frac{1}{Z_{\mrm{U}(N),\Lambda_L,\mbf{B}}}\exp\bigg( S_{\mrm{YM}_{\beta}}(Q)\bigg) dQ,
\end{equs}
where as usual $dQ$ is the product Haar measure, and
\begin{equs}
S_{\mrm{U}(N),\Lambda_L,\mbf{B}}(Q):= N\sum_{p \in P_{\Lambda_L}^+} \beta_p\mrm{Re}\mrm{Tr}(Q_p).
\end{equs}
and $Z_{\mrm{U}(N),\Lambda_L,\mbf{B}}$ is the appropriate normalizing constant.

If $\sup_{p} |\beta_p| < \tilde{\beta}(d)$, then there are constants $C_1=C_1(N,d,|\Lambda_f|,|\Lambda_g|)$, and $C_2=C_2(N,d)$ such that for any smooth local observables $f,g \in C^{\infty}(\mrm{U}(N)^{E_{\Lambda}^+})$, we have that
\begin{equs}
\mathrm{Cov}_{\mu_{\mrm{U}(N),\Lambda_L,\mbf{B}}}(f,g) \leq C_1(\|f\|_{L^{\infty}}+\vertiii{f}_{\infty})(\|g\|_{L^{\infty}}+\vertiii{g}_{\infty}) e^{-C_2 d(\Lambda_f,\Lambda_g)}.
\end{equs}
\end{lemma}

Next, we use Lemma \ref{lm:Mixed-Temp-Mass-Gap} to bound the sensitivity of observables when changing the $\beta_p$ values.

\begin{lemma}[Exponential sensitivity to turning on temperature]\label{lm:temp-gradient-bound}
Suppose that $\sup_p \beta_p < \tilde{\beta}$. Let $f$ be a smooth local observable with $\Lambda_f \sse \Lambda$. Then for any plaquette $p$, we have that
\begin{equs}
    \left|\frac{d}{d\beta_p} \E_{\mrm{YM}_{L,\beta}}[f]\right| \leq C_1(\|f\|_{L^{\infty}}+\vertiii{f}_{\infty})e^{-C_2 d(\Lambda_f,p)}
\end{equs}
for constants $C_1=C_1(N,d,|\Lambda_f|), C_2=C_2(N,d)$. 
\end{lemma}

\begin{proof}
A similar computation as in the proof of Lemma \ref{lm:deriviative-of-conditional-expectation} yields
    \begin{equs}
        \frac{d}{d\beta_p} \E_{\mu_{\mrm{U}(N),\Lambda_L,\mbf{B}}}[f] = \mrm{Cov}_{\mu_{\mrm{U}(N),\Lambda_L,\mbf{B}}}(f(Q),N\mrm{Re}\mrm{Tr}(Q_p)),
    \end{equs}
and the result now follows by Lemma \ref{lm:Mixed-Temp-Mass-Gap}.
\end{proof}

Using Lemma \ref{lm:temp-gradient-bound}, we may now show that the infinite volume limit exists.

\begin{proof}[Proof of the existence infinite volume limit claim in Theorem \ref{thm:U-N-mass-gap}.]
By standard compactness and approximation arguments, it suffices to prove that for any smooth local observable $f$, we have that 
\begin{equ}
\lim_{L \to \infty}\E_{\mu_{\mrm{U}(N),\Lambda_L,\beta}}[f] \text{ exists.}
\end{equ}
Towards this end, let $L_0$ be large enough so that $\Lambda_f \sse \Lambda_{L_0}$. We show that the sequence $\big(\E_{\mu_{\mrm{U}(N),\Lambda_L,\beta}}[f]\big)_{L \geq L_0}$ is Cauchy by establishing the bound:
\begin{equ}
\big|\E_{\mu_{\mrm{U}(N),\Lambda_{L+1},\beta}}[f]-\E_{\mu_{\mrm{U}(N),\Lambda_{L+1},\beta}}[f] \big|\leq C_1 \big(\|f\|_{L^{\infty}}+\vertiii{f}_{\infty} \big)e^{-C_2 d(\Lambda_f,\partial \Lambda_L)}.
\end{equ}
for constants $C_1,C_2$ independent of $L$. The exponential decay implies that the above is summable in $L$, which implies the desired Cauchy property. The idea is to apply the fundamental theorem of calculus and the derivative bound Lemma \ref{lm:temp-gradient-bound}. More precisely, letting $\{p_1,\dots,p_n\}$ be the collection of plaquettes in $\Lambda_{L+1} \backslash \Lambda_L$, and defining $\mbf{B}_k=(\beta^{(k)}_p)_{p \in \mc{P}_{\Lambda_{L+1}^+}}$  by $\beta^{(k)}_p = \beta$ for $p \in P_{\Lambda_L}^{+}$ or $p=p_1,\dots,p_{k-1}$, and $\beta^{(k)}_{p_k}=\beta' \in [0,\beta]$ while $\beta^{(k)}_{p_j}=0$ for $p=p_j$ for $j>k$, we have by Lemma \ref{lm:temp-gradient-bound} that
\begin{equs}
\big|\E_{\mu_{\mrm{U}(N),\Lambda_{L+1},\beta}}[f]-\E_{\mu_{\mrm{U}(N),\Lambda_{L},\beta}}[f] \big| &\leq \sum_{k=1}^{n} \int_{0}^{\beta} \left|\frac{d}{d\beta'} \E_{\mu_{\mrm{U}(N),\Lambda_{L+1},\mbf{B}_k}}[f]\right| d\beta_{p_k}\\
    &\leq n C_1 \big (\|f\|_{L^{\infty}}+\vertiii{f}_{\infty}\big)e^{-C_2 d(\Lambda_f,\partial \Lambda_L)}\\
    &\leq C_1' \big(\|f\|_{L^{\infty}}+\vertiii{f}_{\infty}\big)e^{-C_2' d(\Lambda_f,\partial \Lambda_L)},
\end{equs}
as $n \lesssim L^{d-1}$ while $e^{C_2d(\Lambda_f,\partial \Lambda_L)}$ is exponential in $L$, for $L$ large enough.
\end{proof}

\section{Large $N$ Limit: Proof of Theorem \ref{thm:main-large-N-Limit}}\label{section:Large-N-Limit}

As with the mass gap section, the starting point will be the $\mrm{SU}(N)$ random environment decomposition of Lemma \ref{lm:measure-decompostion}. Applying the decomposition, for any loop $\ell=e_1 e_2 \dots e_n$, we have the following equality in law,
\begin{equs}
    W_{\ell} = e^{i \frac{\theta_{\ell}}{N}} \mrm{tr}(Q_{\ell}),
\end{equs}
where $\theta_{\ell}=\theta_{e_1}+\dots+\theta_{e_{\ell}}$, and $(\theta,Q)$ is distributed according to the measure $\mu_{\mc{U}_L,\beta}$ from definition \ref{def:random-enviornment}. Next recall that the Bakry-Émery criterion verified in Proposition \ref{prop:Poincare-mixing} also implies a Poincare inequality for the measure $\tilde{\mu}_{\mrm{SU}(N),\Lambda_L,\beta}$ (see e.g. \cite[Theorem 5.6.2]{Wang2005Functional}). Namely for any local observable $f$,
\begin{equs}
    \mrm{Var}_{\tilde{\mu}_{\mrm{SU}(N),\Lambda_L,\beta}}(f) \leq \frac{1}{K_{\tilde{S}}}\mc{E}_{\tilde{\mu}_{\mc{U}_L,\beta}}(f,f)=:\sum_{e \in E_{\Lambda}^+} \int \langle \nabla_e f, \nabla_e f\rangle d\tilde{\mu}_{\mrm{SU}(N),\Lambda_L,\beta}.
\end{equs}

In particular setting $f(Q)=\mathrm{tr}(Q_{\ell})$, exactly as computed in \cite[Section 4.2]{shen2023stochastic}, we have the bound $\mc{E}_{\tilde{\mu}_{\mrm{SU}(N),\Lambda_L,\beta}}(f,f) \leq \frac{4n(n-3)}{N}$, so
\begin{equs}
   \mrm{Var}_{_{\tilde{\mu}_{\mrm{SU}(N),\Lambda_L,\beta}}}(f)=\mrm{Var}_{_{\tilde{\mu}_{\mrm{SU}(N),\Lambda_L,\beta}}}(f) \leq \frac{4n(n-3)}{N K_{\tilde{S}}}.
\end{equs}

Next observe that 
\begin{equs}
    |e^{i\frac{\theta_{\ell}}{N}}\mrm{tr}(Q_{\ell})-\mrm{tr}(Q_{\ell})| \leq \sup_{\phi \in [-2\pi n, 2\pi n]}|e^{i \frac{\phi}{N}}-1|.
\end{equs}

Thus,
\begin{equs}
    ~&\mrm{Var}_{\mu_{\mrm{U}(N),\Lambda_L,\beta}}(e^{i\frac{\theta_{\ell}}{N}} \mrm{tr}(Q_{\ell})) \\
    &\leq \mrm{Var}_{\mu_{\mrm{U}(N),\Lambda_L,\beta}}( \mrm{tr}(Q_{\ell}))+\mrm{Var}_{\mu_{\mrm{U}(N),\Lambda_L,\beta}}((e^{i\frac{\theta_{\ell}}{N}} -1)\mrm{tr}(Q_{\ell}))\\
    &\leq \frac{2n(n-3)}{N K_{\tilde{S}}}+ 2\sup_{\phi \in [-2\pi n, 2\pi n]}|e^{i \frac{\phi}{N}}-1|^2
\end{equs}
Finally taking $L \to \infty$ and applying the infinite volume limit part of Theorem \ref{thm:U-N-mass-gap},
\begin{equs}
    \mrm{Var}_{\mu_{\mrm{U}(N),\beta}}(e^{i\frac{\theta_{\ell}}{N}} \mrm{tr}(Q_{\ell})) \leq \frac{4n(n-3)}{N K_{\tilde{S}}}+ 2\sup_{\phi \in [-2\pi n, 2\pi n]}|e^{i \frac{\phi}{N}}-1|^2,
\end{equs}
Now recall that for $\beta<\tilde{\beta}$ , $K_{\tilde{S}} >0$ and is lower bounded by a growing linear function in $N$. So we have $\mrm{Var}_{\mu_{\mrm{U}(N),\beta}}(e^{i\frac{\theta_{\ell}}{N}} \mrm{tr}(Q_{\ell})) \to 0$ as $N \to \infty$, and hence $W_{\ell}-\langle W_{\ell}\rangle_{\mu_{\mrm{U}(N),\beta}} \to 0$ in probability as $N \to \infty$.

Next we prove the factorization property by induction and Cauchy-Schwarz. Fix loops $\ell_1\dots \ell_k$, and let all expectations be taken with respect to $\mu_{\mrm{U}(N),\beta}$.
\begin{equs}
    ~&|\langle W_{\ell_1}W_{\ell_2}\cdots W_{\ell_k}\rangle- \langle W_{\ell_1}\rangle\langle W_{\ell_2}\rangle\cdots \langle W_{\ell_k}\rangle|\\
    &\leq |\langle W_{\ell_1}W_{\ell_2}\cdots W_{\ell_k}\rangle-\langle W_{\ell_1}\rangle\langle W_{\ell_2}\cdots W_{\ell_k}\rangle| \\
    &+ |\langle W_{\ell_1}\rangle\langle W_{\ell_2}\cdots W_{\ell_k}\rangle -\langle W_{\ell_1}\rangle\langle W_{\ell_2}\rangle\cdots \langle W_{\ell_k}\rangle| \\
    & \leq \mrm{Var}(W_{\ell_1})+ |\langle W_{\ell_2}\cdots W_{\ell_k}\rangle -\langle W_{\ell_2}\rangle\cdots \langle W_{\ell_k}\rangle| \to 0
\end{equs}
as $N \to \infty$ where we applied the one loop large $N$ limit result and the inductive hypothesis in the last line, as well as the fact that all Wilson loop observables are bounded by $1$. This completes the proof of Theorem \ref{thm:main-large-N-Limit}.

\bibliographystyle{alpha}
\bibliography{references}

\end{document}